\documentclass[11pt, a4paper]{article}
\usepackage{graphicx}
\usepackage{color}
\usepackage{amsmath}
\usepackage{amssymb}
\usepackage{amscd}
\usepackage{bbm}
\usepackage{tikz}
\usepackage{hyperref}
\hypersetup{
    bookmarks=true,         
    colorlinks=true,       
    linkcolor=red,          
    citecolor=green,        
    filecolor=magenta,      
    urlcolor=cyan           
}
\usepackage[normalem]{ ulem }
\usepackage{soul}
\usepackage{caption}
\usepackage{subcaption}

\newcommand{\R}{\mathbb{R}}
\newcommand{\inr}[1]{\bigl< #1 \bigr>}

\newcommand{\E}{\mathbb{E}}

\newcommand{\eps}{\varepsilon}

\newcommand{\cB}{{\cal B}}
\newcommand{\cM}{{\cal M}}
\newcommand{\cQ}{{\cal Q}}

\newtheorem{Theorem}{Theorem}[section]
\newtheorem{Lemma}[Theorem]{Lemma}

\newtheorem{Definition}[Theorem]{Definition}

\newtheorem{Corollary}[Theorem]{Corollary}
\newtheorem{Remark}[Theorem]{Remark}

\newtheorem{Assumption}{Assumption}[section]

\numberwithin{equation}{section}

\newcommand{\norm}[1]{\left\|#1\right\|}%

\DeclareMathOperator*{\argmin}{argmin}

\def \proof {\noindent {\bf Proof.}\ \ }

\def \endproof
{{\mbox{}\nolinebreak\hfill\rule{2mm}{2mm}\par\medbreak}}

\begin{document}
\title{{Regularization and the small-ball method I: sparse recovery}}
\author{Guillaume Lecu\'e${}^{1,3}$  \and Shahar Mendelson${}^{2,4,5}$}

\footnotetext[1]{CNRS, CREST, ENSAE, 3, avenue Pierre Larousse,
92245 MALAKOFF. France.}
\footnotetext[2]{Department of Mathematics, Technion, I.I.T., Haifa, Israel and Mathematical Sciences Institute, The Australian National University, Canberra, Australia}
 \footnotetext[3] {Email:
guillaume.lecue@ensae.fr }
\footnotetext[4] {Email:
shahar@tx.technion.ac.il}
\footnotetext[5]{Supported by the Israel Science Foundation.}

\maketitle

\begin{abstract}
We obtain bounds on estimation error rates for regularization procedures of the form
\begin{equation*}
  \hat f \in\argmin_{f\in
    F}\left(\frac{1}{N}\sum_{i=1}^N\left(Y_i-f(X_i)\right)^2+\lambda \Psi(f)\right)
\end{equation*}
when $\Psi$ is a norm and $F$ is convex.

Our approach gives a common framework that may be used in the analysis of learning problems and regularization problems alike. In particular, it sheds some light on the role various notions of sparsity have in regularization and on their connection with the size of subdifferentials of $\Psi$ in a neighbourhood of the true minimizer.

As `proof of concept' we extend the known estimates for the LASSO, SLOPE and trace norm regularization.
\end{abstract}

\section{Introduction}\label{sec:intro}
The focus of this article is on {\it regularization}, which is one of the most significant methods in modern statistics. To give some intuition on the method and on the reasons behind its introduction, consider the following standard problem.

Let $(\Omega,\mu)$ be a probability space and set $X$ to be distributed according to $\mu$. $F$ is a class of real-valued functions defined on $\Omega$ and $Y$ is the unknown random variable that one would like to approximate using functions in $F$. Specifically, one would like to identify the best approximation  to $Y$ in $F$, say in the $L_2$ sense, and find the function $f^*$ that minimizes in $F$ the {\it squared loss functional} $f \to \E(f(X)-Y)^2$; that is,
$$
f^*={\rm argmin}_{f \in F} \E(f(X)-Y)^2,
$$
with the underlying assumption that $f^*$ exists and is unique.

Unlike problems in approximation theory, neither the target $Y$ nor the underlying measure $\mu$ are known. Therefore, computing the $L_2$ distance between functions in $F$ and $Y$ is impossible. Instead, one is given partial information: a random sample $(X_i,Y_i)_{i=1}^N$, selected independently according to the joint distribution of $X$ and $Y$.

Because of the random nature of the sample and the limited information it provides, there is no real hope of identifying $f^*$, but rather, only of approximating it. In an {\it estimation problem} one uses the sample to produce a random function $\hat{f} \in F$, and the success of the choice is measured by the distance between $\hat{f}$ and $f^*$ in the $L_2(\mu)$ sense. Thus, one would like to ensure that with high probability with respect to the samples $(X_i,Y_i)_{i=1}^N$, the {\it error rate}
$$
\norm{\hat f-f^*}_{L_2(\mu)}^2 = \E\Big( \big(\hat{f}(X)-f^*(X)\big)^2 | (X_i,Y_i)_{i=1}^N \Big)
$$
is small. More accurately, the question is to identify the way in which the error rate depends on the structure of the class $F$ and scales with the sample size $N$ and the required degree of confidence (probability estimate).

\vskip0.4cm

It is not surprising (and rather straightforward to verify) that the problem becomes harder the larger $F$ is. In contrast, if $F$ is small, chances are that $f^*(X)$ is very far from $Y$, and identifying it, let alone approximating it, is pointless.

In situations we shall refer to as {\it learning problems}, the underlying assumption is that $F$ is indeed small, and the issue of the approximation error -- the distance between $Y$ and $f^*$ is ignored.

While the analysis of learning problems is an important and well-studied topic, the assumption that $F$ is reasonably small seems somewhat restrictive; it certainly does not eliminate the need for methods that allow one to deal with very large classes.

Regularization was introduced as an alternative to the assumption on the `size' of $F$. One may consider large classes, but combine it with the belief that $f^*$ belongs to a relatively small substructure in $F$. The idea is to penalize a choice of a function that is far from that substructure, which forces the learner to choose a function in the `right part' of $F$.

\vskip0.4cm

Formally, let $E$ be a vector space, assume that $F \subset E$ is a closed and convex set and let $\Psi:E \to \R^+$ be the penalty. Here, we will only consider the case in which $\Psi$ is a norm on $E$.

Let $\lambda>0$ and for a sample $(X_i,Y_i)_{i=1}^N$, set
\begin{equation*}
  \hat f \in\argmin_{f\in
    F}\left(\frac{1}{N}\sum_{i=1}^N\left(Y_i-f(X_i)\right)^2+\lambda \Psi(f)\right);
\end{equation*}
$\hat{f}$ is called a regularization procedure, $\Psi$ is the regularization function and $\lambda$ is the regularization parameter.

\vskip0,4cm

In the classical approach to regularization, the substructure of $f^*$ is quantified directly by $\Psi$. The underlying belief is that $\Psi(f^*)$ is not `too big' and one expects the procedure to produce $\hat{f}$ for which $\Psi(\hat{f})$ is of the order of $\Psi(f^*)$. Moreover, the anticipated error rate $\|\hat{f}-f^*\|_{L_2(\mu)}$ depends on $\Psi(f^*)$. In fact, an optimistic viewpoint is that regularization could perform as well as the best learning procedure in the class $\{f : \Psi(f) \leq \Psi(f^*)\}$, but without knowing $\Psi(f^*)$ beforehand.

Among the regularization schemes that are based on the classical approach are reproducing kernel Hilbert spaces (RKHS), in which the RKHS norm serves as the penalty. Since RKHS norms capture various notions of smoothness, in RKHS regularization one is driven towards a choice of a smooth $\hat{f}$ -- as smooth as $f^*$ is.

\vskip0.4cm

In more modern regularization problems the situation is very different. Even when penalizing with a norm $\Psi$, one no longer cares whether or not $\Psi(f^*)$ is small; rather, one knows (or at least believes) that $f^*$ is {\it sparse} in some sense, and the hope is that this sparsity will be reflected in the error rate.

In other words, although one uses certain norms as regularization functions -- norms that seemingly have nothing to do with `sparsity' -- the hope is that the sparse nature of $f^*$ will be exposed by the regularization procedure, while $\Psi(f^*)$ will be of little importance.

\vskip0.4cm

The most significant example in the context of sparsity-driven regularization is the celebrated LASSO estimator \cite{MR1379242}.
Let $F = \{\inr{t,\cdot} :t \in \R^d\}$ and set $t^*$ to be a minimizer in $\R^d$ of the functional $t \to \E(\inr{t,X}-Y)^2$. The LASSO is defined by
$$
\hat{t}\in\argmin_{t \in \R^d} \Big(\frac{1}{N}\sum_{i=1}^N \left(\inr{t,X_i}-Y_i\right)^2 + \lambda \Psi(t) \Big)
$$
for the choice  $\Psi(t)=\|t\|_1 =\sum_{i=1}^d |t_i|$.

The remarkable property of the LASSO (see \cite{MR2829871} and \cite{MR2533469}) is that for a well-chosen regularization parameter $\lambda$, if $t^*$ is supported on at most $s$ coordinates (and under various assumptions on $X$ and $Y$ to which we will return later), then with high probability,
$$
\|\hat{t}-t^*\|_2^2 \lesssim  \frac{s \log(ed)}{N}.
$$
Thus, the error rate of the LASSO does not depend on $\Psi(t^*)=\|t^*\|_1$, but rather on the degree of sparsity of $t^*$, measured here by the cardinality of its support $\norm{t^*}_0=|\{i : t_i^* \not = 0\}|$.

This fact seems almost magical, because to the naked eye, the regularization function $\|t\|_1$ has nothing to do with sparsity; yet $\ell_1$ regularization leads to a sparsity-driven error rate.

A standard (yet somewhat unconvincing) explanation of this phenomenon is that the penalty $\|t\|_1$ is a convexified version of $\|t\|_0=|\{i : t_i \not = 0\}|$, though this loose connection hardly explains why $\|t^*\|_0$ has any effect on the error rate of the LASSO.

A similar phenomenon occurs for other choices of $\Psi$, such as the SLOPE and trace-norm regularization, which will be explored in detail in what follows. In all these cases and others like them, the regularization function is a norm that does not appear to be connected to sparsity, nor to other natural notions of low-dimensional structures for that matter. Yet, and quite mysteriously, the respective regularization procedure emphasizes those very properties of $t^*$.

The aim of this note is to offer a framework that can be used to tackle standard learning problems (small $F$) and regularized problems alike. Moreover, using the framework, one may explain how certain norms lead to the emergence of sparsity-based bounds.

\vskip0.4cm

In what follows we will show that two parameters determine the error rate of regularization problems. The first one captures the `complexity' of each set in the natural hierarchy in $F$
$$
F_\rho=\{f \in F : \Psi(f-f^*) \leq \rho\}.
$$
Applying results from \cite{Shahar-COLT,Shahar-ACM,shahar_general_loss}, the `complexity' of each $F_\rho$ turns out to be the optimal (in the minimax sense) error rate of the learning problem in that set. To be more precise, the main ingredient in obtaining a sharp error rate of a learning problem in a class $H$ is an accurate analysis of the empirical excess squared loss functional
\begin{equation}\label{eq:excess-loss}
f \to P_N {\cal L}_f = \frac{1}{N}\sum_{i=1}^N (f(X_i)-Y_i)^2 - \frac{1}{N}\sum_{i=1}^N (f^*(X_i)-Y_i)^2.
\end{equation}
Since the minimizer $\hat{f}$ of the functional \eqref{eq:excess-loss} satisfies $P_N {\cal L}_{\hat{f}} \leq 0$, one may obtain an estimate on the error rate by showing that with high probability, if $\|f-f^*\|_{L_2(\mu)} \geq r$ then $P_N {\cal L}_f >0$. This excludes functions in the set $\{f \in H: \|f-f^*\|_{L_2(\mu)} \geq r\}$ as potential empirical minimizers.
That `critical level' turns out to be the correct (minimax) error rate of a learning problem in $H$. That very same parameter is of central importance in regularization problems --- specifically, the `critical level' $r(\rho)$ for each one of the sets $\{f \in F : \Psi(f-f^*) \leq \rho\}$ (see Section \ref{sec:decomp_excess_loss} for an accurate definition of $r(\rho)$ and its role in the analysis of learning problems and regularization problems).

\vskip0.4cm

The second parameter, which is the main ingredient in our analysis of regularization problems, measures the `size' of the subdifferential of $\Psi$ in points that are close to $f^*$: recall that the subdifferential of $\Psi$ in $f$ is
\begin{equation*}
(\partial\Psi)_f = \{z^*\in E^*: \Psi(f+h)\geq \Psi(f)+z^*(h) \mbox{ for every } h\in E  \}
\end{equation*}
where $E^*$ is the dual space of the normed space $(E,\Psi)$, and that if $f\neq0$, the subdifferential consists of all the norm one linear functionals $z^*$ for which $z^*(f)=\Psi(f)$.

Fix $\rho>0$ and let $\Gamma_{f^*}(\rho)$ be the collection of functionals that belong to the subdifferential $(\partial \Psi)_f $ for some $f \in F$ that satisfies $\Psi(f-f^*) \leq \rho/20$.
Set
$$
H_\rho = \{f \in F : \Psi(f-f^*)=\rho \ {\rm and} \ \|f-f^*\|_{L_2(\mu)} \leq r(\rho)\}
$$
and let
\begin{equation*}
\Delta(\rho) = \inf_{h \in H_\rho} \sup_{z^* \in \Gamma_{f^*}(\rho)} z^*(h-f^*).
\end{equation*}
Hence, $\Gamma_{f^*}(\rho)$ is a subset of the unit sphere of $E^*$ when $0\notin \{f\in F:\Psi(f-f^*)\leq \rho/20\}$ and it is the entire unit ball of $E^*$ otherwise. And, since $H_\rho$ consists of functions whose $\Psi$ norm is $\rho$, it is evident that $\Delta(\rho) \leq \rho$. Therefore, if $\Delta(\rho) \geq \alpha \rho$ for a fixed $0<\alpha \leq 1$ then $\Gamma_{f^*}(\rho)$ is rather large: for every $h \in H_\rho$ there is some $z^* \in \Gamma_{f^*}(\rho)$ for which $z^*(h)$ is `almost extremal'---that is, at least $\alpha \rho$.

\vskip0.4cm

{\bf Our main result (Theorem \ref{lemma:basic-combining-loss-and-reg} below) is that if $\Gamma_{f^*}(\rho)$ is large enough to ensure that $\Delta(\rho) \geq 4\rho/5$, and the regularization parameter $\lambda$ is set to be of the order of $\frac{r^2(\rho)}{\rho}$, then with high probability, the regularized minimizer in $F$, $\hat{f}$, satisfies that $\|\hat{f}-f^*\|_{L_2(\mu)} \leq r(\rho)$ and $\Psi(\hat f-f^*)\leq \rho$.}

\vskip0.4cm

Theorem \ref{lemma:basic-combining-loss-and-reg} implies that one may analyze regularization problems by selecting $\rho$ wisely, keeping in mind that points in a $\Psi$-ball of radius $\sim \rho$ around $f^*$ must generate a sufficiently large subdifferential. And the fact that functionals in $\Gamma_{f^*}(\rho)$ need to be `almost extremal' only for points in $H_\rho$ rather than for the entire sphere is crucial; otherwise, it would have forced $\Gamma_{f^*}(\rho)$ to be unreasonably large -- close to the entire dual sphere.


\vskip0.4cm

As will be clarified in what follow, sparsity, combined with the right choice of $\Psi$, contributes in two places: firstly, if $f^*$ is sparse in some sense and $\Psi$ is not smooth on sparse elements, then $\Gamma_{f^*}(\rho)$, which contains the subdifferential $(\partial \Psi)_{f^*}$, is large; secondly, for the right choice of $\rho$ the `localization' $H_\rho$ consists of elements that are well placed: if $\Psi(f-f^*)=\rho$ and $\|f-f^*\|_{L_2(\mu)} \leq r(\rho)$, there is some $z^* \in \Gamma_{f^*}(\rho)$ for which $z^*(f-f^*)$ is large enough. The fact that $H_\rho$ is well placed is an outcome of some compatibility between $\Psi$ and the $L_2(\mu)$ norm.

Of course, to find the right choice of $\rho$ one must first identify $r(\rho)$, which is, in itself, a well-studied yet nontrivial problem.

\vskip0.4cm

Before we dive into technical details, let us formulate some outcomes of our main result. We will show how it can be used to obtain sparsity-driven error rates in three regularization procedures: the LASSO, SLOPE and trace norm regularization. In all three cases our results actually extend the known estimates in various directions.

\vskip0.4cm

\noindent{\bf The LASSO.}

The LASSO is defined for the class of linear functional $F=\{\inr{t,\cdot} : t \in \R^d\}$. For a fixed $t_0 \in \R^d$, the goal is to identify $t_0$ using linear measurements, the regularization function is $\Psi(t)=\|t\|_1 = \sum_{i=1}^d |t_i|,$
and the resulting regularization procedure produces
$$
\hat{t} \in\argmin_{t \in \R^d} \Big(\frac{1}{N}\sum_{i=1}^N \left(\inr{t,X_i}-Y_i\right)^2 + \lambda \|t\|_1 \Big).
$$
The LASSO has been studied extensively in the last two decades. Even though some recent advances \cite{MR3153940,MR3224285,MR3161450} have shown the LASSO to have its limitation, historically, it has been the benchmark estimator of high-dimensional statistics --- mainly because a high dimensional parameter space does not significantly affect its performance as long as $t_0$ is sparse. This was shown for example, in \cite{MR2533469,MR2386087,deter_lasso,MR2396809,MR2488351,MR3025133,MR3181133} in the context of estimation and sparse oracle inequalities, in \cite{MR2278363,MR2274449,bach_support} for support recovery results;  and in various other instances as well; we refer the reader to the books \cite{MR2807761,MR2829871} for more results and references on the LASSO.

\vskip0.4cm

\noindent{\bf SLOPE.}

In some sense, SLOPE, introduced in \cite{slope1,slope2}, is actually an extension of the LASSO, even though it has been introduced as an extension of multiple-test procedures. Again, the underlying class is $F=\{\inr{t,\cdot} : t \in \R^d\}$, and to define the regularization function
let $\beta_1 \geq \beta_2 \geq ... \geq \beta_d>0$ and set
$$
\Psi(t)=\sum_{i=1}^d \beta_i t_i^\sharp,
$$
where $(t_i^\sharp)_{i=1}^d$ denotes the non-increasing re-arrangement of $(|t_i|)_{i=1}^d$. Thus, the SLOPE norm is a sorted, weighted $\ell_1$-norm, and for $(\beta_1,...,\beta_d)=(1,...,1)$, SLOPE regularization coincides with the LASSO.

\vskip0.4cm

\noindent{\bf Trace-norm regularization.}

Consider the trace inner-product on $\R^{m \times T}$. Let $F=\{\inr{A,\cdot} : A \in \R^{m \times T}\}$ and given a target $Y$ put $A^*$ to be the matrix that minimizes $A \to\E (\inr{A,X}-Y)^2$. The regularization function is the {\it trace norm}.
\begin{Definition}
Let $A$ be a matrix and set $(\sigma_i(A))$ to be its singular values, arranged in a non-increasing order.
For $p \geq 1$, $\|A\|_p = (\sum \sigma_i^p(A))^{1/p}$ is the $p$-Schatten norm. 
\end{Definition}
Note that the trace-norm is simply the $1$-Schatten norm, the Hilbert-Schmidt norm is the $2$-Schatten norm and the operator norm is the $\infty$-Schatten norm.

The trace norm regularization procedure is
\begin{equation*}
 \hat A \in\argmin_{A\in\R^{m\times T}}\Big(\frac{1}{N}\sum_{i=1}^N (Y_i-\inr{X_i,A})^2+\lambda \|A\|_1\Big)
 \end{equation*}
and it was introduced for the reconstruction of low-rank, high-dimensional matrices  \cite{MR2680543,MR2906869,MR2816342,MR2809094,MR2815834,MR2930649}.

\vskip0.5cm

As will be explained in what follows, our main result holds in rather general situations and may be implemented in examples once the `critical levels' $r(\rho)$ are identified. Since the examples we present serve mainly as ``proof of concept", we will focus only on one scenario in which $r(\rho)$ may be completely characterized for an arbitrary class of functions.

\begin{Definition} \label{def:subgaussian}
Let $\ell_2^M$ be an $M$-dimensional inner product space and let $\mu$ be a measure on $\ell_2^M$.
The measure $\mu$ is isotropic if for every $x \in \ell_2^M$,
$$
\int \inr{x,t}^2 d\mu(t) = \|x\|_{\ell_2^M}^2;
$$
it is $L$-subgaussian if for every $p \geq 2$ and every $x \in \ell_2^M$,
$$
\|\inr{x,\cdot}\|_{L_p(\mu)} \leq L\sqrt{p} \|\inr{x,\cdot}\|_{L_2(\mu)}.
$$
\end{Definition}
Hence, the covariance structure of an isotropic measure coincides with the inner product in $\ell_2^M$, and if $\mu$ is an $L$-subgaussian measure then the $L_p(\mu)$ norm of a linear form  does not grow faster than the $L_p$ norm of the corresponding Gaussian variable.

\begin{Assumption}\label{ass:examples}
Assume that the underlying measure $\mu$ is isotropic and $L$-subgaussian, and that for $f^*=\inr{t^*,\cdot}$ (or $f^*=\inr{A^*,\cdot}$ in the matrix case), the noise\footnote{In what follows we will refer to $\xi$ as `the noise' even though it depends in general on $Y$ and $X$. The reason for using that term comes from the situation in which $Y=f^*(X)-W$ for a symmetric random variable $W$ that is independent of $X$ (independent additive noise); thus $\xi=W$. We have opted to call $\xi$ `the noise' because its role in the general case and its impact on the error rate is rather similar to what happens for independent noise.} $\xi=f^*(X)-Y$ belongs to $L_q$ for some $q>2$.
\end{Assumption}

When dealing with the LASSO and SLOPE, the natural Euclidean structure is the standard one in $\R^d$, and for trace norm regularization, the natural Euclidean structure is endowed by the trace inner product in $\R^{m \times T}$.

\begin{Remark}\label{rem:subgauss_isotropicicty}
In the supplementary material we study a general $X$ without assuming it is isotropic, which means dealing with less natural Euclidean structures in the examples we present. It is also possible to go beyond the subgaussian case, we refer the reader to \cite{LM_reg_comp} where other moment assumptions on $X$ are considered.
\end{Remark}

\vskip0.4cm

The second part of Assumption \ref{ass:examples}, that $\xi \in L_q$ for some $q>2$, is rather minimal. Indeed, for the functional $f \to \E (f(X)-Y)^2$ to be well defined, one must assume that $f(X)-Y \in L_2$; the assumption here is only slightly stronger.

\vskip0.4cm

Applying our main result we will show the following:
\begin{Theorem} \label{thm:intro-LASSO-est}
Consider the LASSO under Assumption \ref{ass:examples}. Let $0<\delta<1$. Assume that there is some $v \in \R^d$ supported on at most $s$ coordinates for which
	\begin{equation*}
	 \norm{t^*-v}_1\leq c_1(\delta) \|\xi\|_{L_q} s \sqrt{\frac{\log(ed)}{N}}.
	 \end{equation*}
If $\lambda= c_2(L,\delta)\|\xi\|_{L_q} \sqrt{\log(ed)/N}$ and $N\geq s \log(ed/s)$, then with probability at least $1-\delta$ the LASSO estimator with regularization parameter $\lambda$ satisfies that
for every $1\leq p \leq 2$
\begin{equation*}
\norm{\hat t-t^*}_p \leq c_3(L,\delta)\|\xi\|_{L_q} s^{1/p}\sqrt{\frac{\log(ed)}{N}}.
\end{equation*}
\end{Theorem}
The error rate in Theorem \ref{thm:intro-LASSO-est} coincides with the standard estimate on the LASSO (cf. \cite{MR2533469}), but in a broader context: $t^*$ need not be sparse but only approximated by a sparse vector; the target $Y$ is arbitrary and the noise $\xi$ may be heavy tailed and need not be independent of $X$.

\vskip0.4cm
Turning to SLOPE, let us recall the estimates from \cite{slope2}, where the setup is somewhat restricted: Let $X$ be a Gaussian vector on $\R^d$, set $W$ to be a Gaussian random variable with variance $\sigma^2$ that is independent of $X$ and put $Y=\inr{t^*,X}+W$.
Consider some $q\in(0,1)$, let $\Phi^{-1}(\alpha)$ be the $\alpha$-th quantile of the standard normal distribution and put $\beta_i=\Phi^{-1}(1-iq/(2d))$.

\begin{Theorem} \cite{slope2} \label{theo:candes_slope}
Let $1\leq s\leq d$ satisfy that $s/d=o(1)$ and $(s \log d)/N=o(1)$ when $N\rightarrow \infty$. If $0<\eps<1$, $N \rightarrow \infty$ and $\lambda = 2\sigma/\sqrt{N}$, the SLOPE estimator with weights $(\beta_i)_{i=1}^d$ and regularization parameter $\lambda$ satisfies
\begin{equation*}
  \sup_{\norm{t^*}_0\leq s}Pr\Big(\frac{N\norm{\hat
      t-t^*}_2^2}{2\sigma^2 s \log(d/s)}>1+3\eps\Big)\rightarrow 0.
\end{equation*}
\end{Theorem}

Note that Theorem \ref{theo:candes_slope} is asymptotic in nature and not `high-dimensional'. Moreover, it only holds for a Gaussian $X$, independent Gaussian noise $W$, a specific choice of weights $(\beta_i)_{i=1}^d$ and $t^*$ that is $s$-sparse.

\vskip0.4cm

We consider a more general situation. Let $\beta_i \leq C \sqrt{\log(ed/i)}$ and set $\Psi(t)=\sum_{i=1}^d  \beta_i t_i^\sharp $.

\begin{Theorem} \label{thm:intro-SLOPE-est}
There exists constants $c_1,c_2$ and $c_3$ that depend only on $L$, $\delta$ and $C$ for which the following holds. Under Assumption \ref{ass:examples}, if there is $v\in\R^d$ that satisfies $|{\rm supp}(v)|\leq s$ and
	\begin{equation*}
	 \Psi(t^*-v)\leq c_1 \|\xi\|_{L_q} \frac{ s}{\sqrt{N}} \log\Big(\frac{ed}{s}\Big),
	 \end{equation*}
then for $N \geq  c_2 s \log(ed/s)$ and with the choice of $\lambda= c_2 \|\xi\|_{L_q}/\sqrt{N}$, one has
\begin{equation*}
\Psi(\hat t-t^*)\leq c_3\|\xi\|_{L_q} \frac{s}{\sqrt{N}} \log\Big(\frac{ed}{s}\Big) \ \ \mbox{ and } \ \ \norm{\hat t-t^*}_2^2 \leq c_3 \|\xi\|_{L_q}^2 \frac{s }{N}\log\Big(\frac{ed}{s}\Big)
\end{equation*}
with probability at least $1-\delta$.
\end{Theorem}

Finally, let us consider trace norm regularization.
\begin{Theorem} \label{thm:intro-TRACE-est}
Under Assumption \ref{ass:examples} and if there is $V\in\R^{m\times T}$ that satisfies that ${\rm rank}(V)\leq s$ and
	\begin{equation*}
	 \norm{A^*-V}_1\leq c_1\|\xi\|_{L_q} s \sqrt{\frac{\max\{m,T\}}{N}},
	 \end{equation*}
one has the following.
Let $N \geq c_2 s \max\{m,T\}$ and $\lambda = c_3\|\xi\|_{L_q} \sqrt{\frac{\max\{m,T\}}{N}}$. Then with probability at least $1-\delta$, for any $1\leq p\leq2$
\begin{equation*}
\norm{\hat A-A^*}_p \leq  c_4\|\xi\|_{L_q}  s^{1/p}\sqrt{\frac{\max\{m,T\}}{N}}.
\end{equation*}
The constants $c_1,c_2,c_3$ and $c_4$ depends only on $L$ and $\delta$.
\end{Theorem}

A result of a similar flavour to Theorem \ref{thm:intro-TRACE-est} is Theorem~9.2 from \cite{MR2829871}.
\begin{Theorem} \label{thm:Kolt-intro}
Let $X$ be an isotropic and $L$-subgaussian vector, and $W$ that is mean-zero, independent of $X$ and belongs to the Orlicz space $L_{\psi_\alpha}$ for some $\alpha \geq 1$. If $Y=\inr{A^*,X}+W$ and
\begin{equation*}
  \lambda\geq c_1(L) \max\left\{\norm{\xi}_2\sqrt{\frac{m(t+\log m)}{N}},  \norm{\xi}_{\psi_\alpha}\log^{1/\alpha}\Big(\frac{\norm{\xi}_{\psi_\alpha}}{\norm{\xi}_{L_2}}\Big)\frac{\sqrt{m}(t
    + \log N)(t+\log m)}{N}\right\},
\end{equation*}
then with probability at least $1-3\exp(-t)-\exp(-c_2(L) N)$
\begin{equation}
  \label{eq:vlad}
  \norm{\hat A-A^*}_{2}^2\leq c_3 \min\left\{\lambda
\norm{A^*}_{1}, \lambda^2{\rm rank }(A^*)\right\}.
\end{equation}
\end{Theorem}

Clearly, the assumptions of Theorem \ref{thm:Kolt-intro} are more restrictive than those of Theorem \ref{thm:intro-TRACE-est}, as the latter holds for a heavy tailed $\xi$ that need not be independent of $X$, and for $A^*$ that can be approximated by a low-rank matrix. Moreover, if $\|A^*\|_1$ is relatively large and the error rate in Theorem \ref{thm:Kolt-intro} is the sparsity-dominated $\lambda^2 {\rm rank}(A^*)$, then the error rate in Theorem \ref{thm:intro-TRACE-est} is better by a logarithmic factor.

\vskip0.4cm

The proofs of the error rates in all the three examples will be presented in Section \ref{sec:r}.

\subsection{Notation}
We end the introduction with some standard notation.

Throughout, absolute constants are denoted by $c,c_1...$, etc. Their value may change from line to line. When a constant depends on a parameter $\alpha$ it will be denoted by $c(\alpha)$. $A \lesssim B$ means that $A \leq cB$ for an absolute constant $c$, and the analogous two-sided inequality is denoted by $A \sim B$. In a similar fashion, $A \lesssim_\alpha B$ implies that $A \leq c(\alpha)B$, etc.

Let $E \subset L_2(\mu)$ be a vector space and set $\Psi$ to be a norm on $E$. For a set $A \subset E$, $t \in E$ and $r>0$, let $rA+t=\{ra+t : a \in A\}$.

Denote by $B_\Psi=\{w \in E: \Psi(w) \leq 1\}$ the unit ball of $(E,\Psi)$ and set $S_\Psi=\{f\in E : \Psi(f) =1\}$ to be the corresponding unit sphere. $B_\Psi(\rho,f)$ is the ball of radius $\rho$ centred in $f$ and $S_\Psi(\rho,f)$ is the corresponding sphere. Also, set $D$ to be the unit ball in $L_2(\mu)$, $S$ is the unit sphere there, and $D(\rho,f)$ and $S(\rho,f)$ are the ball and sphere centred in $f$ and of radius $\rho$, respectively.

A class of spaces we will be interested in consists of $\ell_p^d$, that is, $\R^d$ endowed with the $\ell_p$ norm; $B_p^d$ denotes the unit ball in $\ell_p^d$ and  $S(\ell_p^d)$ is the unit sphere.

For every $x=(x_i)_{i=1}^d$,  $(x_i^\sharp)_{i=1}^d$ denotes the non-increasing rearrangement of $(|x_i|)_{i=1}^d$.

Finally, if $(X_i,Y_i)_{i=1}^N$ is a sample,
$P_N h = \frac{1}{N}\sum_{i=1}^N h(X_i,Y_i)$
is the empirical mean of $h$.

\section{Preliminaries: The regularized functional}\label{sec:reg_func}

Let $F \subset E$ be a closed and convex class of functions. Recall that for target $Y$, $f^*$ is the minimizer in $F$ of the functional $f \to \E(f(X)-Y)^2$. Since $F$ is closed and convex, the minimum exists and is unique.

Let ${\cal L}_f(X,Y) = (f(X)-Y)^2-(f^*(X)-Y)^2$ be the excess squared loss functional and for $\lambda>0$ let
$$
{\cal L}_f^\lambda(X,Y) = {\cal L}_f + \lambda(\Psi(f)-\Psi(f^*))
$$
be its regularized counterpart. Thus, for a random sample $(X_i,Y_i)_{i=1}^N$, the empirical (regularized) excess loss functional is
$$
P_N {\cal L}_f^\lambda = \frac{1}{N}\sum_{i=1}^N {\cal L}_f(X_i,Y_i) + \lambda(\Psi(f)-\Psi(f^*)),
$$

Note that if $\ell_f(x,y)=(y-f(x))^2$ and $\hat{f}$ minimizes $P_N \ell_f + \lambda \Psi(f)$ then $\hat{f}$ also minimizes $P_N{\cal L}_f^\lambda$. Moreover, since ${\cal L}^{\lambda}_{f^*}=0$, it is evident that $P_N {\cal L}_{\hat{f}}^{\lambda} \leq 0$.

This simple observation shows that the random set $\{f \in F : P_N {\cal L}_f^\lambda > 0\}$ may be excluded from our considerations, as it does not contain potential minimizers. Therefore, if one can show that with high probability,
$$
\{f \in F : P_N {\cal L}_f^\lambda \leq  0\} \subset \{f \in F : \|f-f^*\|_{L_2(\mu)} \leq r\},
$$
then on that event, $\|\hat{f}-f^*\|_{L_2(\mu)} \leq r$.

We will identify when $P_N {\cal L}_f^\lambda>0$ by considering the two parts of the empirical functional:
the empirical excess loss $P_N {\cal L}_f$ and the regularized part $\lambda(\Psi(f)-\Psi(f^*))$.

\vskip0.4cm

Because of its crucial role in obtaining error estimates in learning problems, the functional $f \to P_N {\cal L}_f$ has been studied extensively using the {\it small-ball method}, (see, e.g., \cite{Shahar-COLT,Shahar-ACM,shahar_general_loss}). Thus, the first component in the machinery we require for explaining both learning problems and regularization problems is well understood and ready-to-use; its details are outlined below.

\vskip0.4cm
\subsection{The natural decomposition of $P_N{\cal L}_f$}\label{sec:decomp_excess_loss}
Set $\xi = \xi(X,Y)=f^*(X)-Y$ and observe that
\begin{align*}
{\cal L}_f(X,Y)= & (f-f^*)^2(X)+2(f-f^*)(X) \cdot (f^*(X)-Y)
\\
= & (f-f^*)^2(X) + 2\xi(f-f^*)(X).
\end{align*}
Since $F$ is convex, the characterization of the nearest point map in a Hilbert space shows that
$$
\E (f-f^*)(X) \cdot (f^*(X)-Y) \geq 0
$$
 for every $f \in F$. Hence, setting $\xi_i=f^*(X_i)-Y_i$, one has
\begin{align*}
P_N {\cal L}_f^\lambda \geq & \frac{1}{N} \sum_{i=1}^N (f-f^*)^2(X_i) +  2\Big(\frac{1}{N} \sum_{i=1}^N \xi_i (f-f^*)(X_i) - \E \xi (f-f^*)(X)\Big)
\\
+ & \lambda(\Psi(f)-\Psi(f^*)).
\end{align*}
To simplify notation, for $w \in L_2(\mu)$ set
${\cal Q}_{w} = w^2$ and ${\cal M}_{w} = \xi w - \E \xi w$.
Thus, for every $f \in F$,
\begin{equation} \label{eq:quad-multi-decomposition}
P_N {\cal L}_f^\lambda \geq P_N {\cal Q}_{f-f^*} +2 P_N {\cal M}_{f-f^*} + \lambda(\Psi(f)-\Psi(f^*)).
\end{equation}

The decomposition of the empirical excess loss to the quadratic component (${\cal Q}_{f-f^*}$) and the multiplier one (${\cal M}_{f-f^*}$) is the first step in applying the small-ball method to learning problems. One may show that on a large event, if $\|f-f^*\|_{L_2(\mu)}$ is larger than some critical level then $P_N {\cal Q}_{f-f^*} \geq \theta \|f-f^*\|_{L_2}^2$ and dominates $P_N {\cal M}_{f-f^*}$; hence
$P_N {\cal L}_f > 0$.
\vskip0.4cm

To identify this critical level, let us define the following parameters:
\begin{Definition}\label{def:the-three-paraemters}
Let $H \subset F$ be a convex class that contains $f^*$. Let $(\eps_i)_{i=1}^N$ be independent, symmetric, $\{-1,1\}$-valued random variables that are independent of $(X_i,Y_i)_{i=1}^N$.

For $\gamma_Q,\gamma_M>0$ set
\begin{equation*}
r_Q(H,\gamma_Q) = \inf\left\{r>0: \E \sup_{h \in H \cap D(r,f^*)} \left|\frac{1}{N}\sum_{i=1}^N \eps_i (h-f^*)(X_i)\right| \leq \gamma_Q r \right\},
\end{equation*}
let
\begin{equation*}
\phi_N(H,s)= \sup_{h \in H \cap D(s,f^*)} \left|\frac{1}{\sqrt{N}} \sum_{i=1}^N \eps_i \xi_i (h-f^*)(X_i)\right|,
\end{equation*}
and put
\begin{equation*}
r_M(H,\gamma_M,\delta) =  \inf\left\{ s>0 : Pr \left(  \phi_N(H,s)\leq \gamma_M s^2 \sqrt{N} \right) \geq 1-\delta \right\}.
\end{equation*}
\end{Definition}
The main outcome of the small-ball method is that for the right choices of $\gamma_M$ and $\gamma_Q$, $r=\max\{r_M,r_Q\}$ is the above-mentioned `critical level' in $H$, once $H$ satisfies a weak small-ball condition.

\begin{Assumption}[The small ball condition]\label{ass:small-ball}
Assume that there are constants $\kappa>0$ and $0<\eps \leq 1$, for which, for every $f,h \in F \cup \{0\}$,
$$
Pr\left(|f-h| \geq \kappa \|f-h\|_{L_2(\mu)}\right) \geq \eps.
$$
\end{Assumption}
There are numerous examples in which the small-ball condition may be
verified for $\kappa$ and $\eps$ that are absolute constants. We refer the reader to
\cite{LM_compressed,shahar_general_loss,Shahar-Vladimir,MR3364699,Shahar-ACM,RV_small_ball}
for some of them.

\begin{Theorem}[\cite{Shahar-ACM}]\label{thm:small-ball-method}
Let $H$ be a closed, convex class of functions that contains $f^*$ and satisfies Assumption \ref{ass:small-ball} with constants $\kappa$ and $\eps$. If
$\theta = \kappa^2 \eps/16$ then for every $0<\delta<1$, with probability at least $1-\delta-2\exp(-N\eps^2/2)$ one has:
\begin{description}
\item{$\bullet$} for every $f \in H$,
$$
|P_N {\cal M}_{f-f^*}| \leq \frac{\theta}{8}\max\left\{\|f-f^*\|_{L_2(\mu)}^2,r_M^2 \left(H,\theta/10,\delta/4\right)\right\},
$$
\item{$\bullet$} If $f \in H$ and $\|f-f^*\|_{L_2(\mu)} \geq r_Q \left(H,\kappa \eps/32\right)$ then
$$
P_N {\cal Q}_{f-f^*} \geq \theta \|f-f^*\|_{L_2(\mu)}^2.
$$
\end{description}
In particular, with probability at least $1-\delta-2\exp(-N\eps^2/2)$,
$$
P_N {\cal L}_f \geq \frac{\theta}{2}
  \|f-f^*\|_{L_2(\mu)}^2
$$
for every $f\in H$ that satisfies
$$
\|f-f^*\|_{L_2(\mu)} \geq \max\left\{r_M \left(H,\theta/10,\delta/4\right),r_Q \left(H,\kappa \eps/32\right)\right\}.
$$
\end{Theorem}

From now on, we will assume that $F$ satisfies the small-ball condition with constants $\kappa$ and $\eps$, and that $\theta=\kappa^2 \eps/16$.
\begin{Definition}
Let $\rho>0$ and set
\begin{equation*}
r_M(\rho)=  r_M \bigl(F \cap B_\Psi(\rho,f^*) ,\frac{\theta}{10},\frac{\delta}{4}\bigr) \ \
\mbox{ and } \ \
r_Q(\rho) =  r_Q \bigl(F \cap B_\Psi(\rho,f^*),\frac{\kappa \eps}{32}\bigr).
\end{equation*}
In what follows we will abuse notation and omit the dependence of $r_M$ and $r_Q$ on $f^*$, $\kappa$, $\eps$ and $\delta$.

Let $r(\cdot)$ be a function that satisfies $r(\rho)\geq \sup_{f^* \in F} \max\{r_Q(\rho),r_M(\rho)\}.$ Finally, put
\begin{equation*}
{\cal O}(\rho) = \sup_{f\in F \cap B_\Psi(\rho, f^*)\cap D(r(\rho), f^*)}\big| P_N \cM_{f-f^*}\big|.
\end{equation*}
\end{Definition}
Theorem \ref{thm:small-ball-method} implies the following:
\begin{Corollary}[\cite{Shahar-ACM}]\label{cor:small-ball-with-O}
Using the notation introduced above, on an event of probability at least $1-\delta-2\exp(-N\eps^2/2)$, if $f \in F \cap B_\Psi(\rho,f^*)$ and $\|f-f^*\|_{L_2(\mu)} \geq r(\rho)$ then
$$
P_N {\cal L}_f \geq \frac{\theta}{2}\|f-f^*\|_{L_2(\mu)}^2.
$$
Moreover, on the same event,
$$
{\cal O}(\rho)\leq \frac{\theta}{8}r^2(\rho).
$$
\end{Corollary}

\begin{Remark}
Let us stress once again that $r(\rho)$  plays a central role in the analysis of {\it empirical risk minimization} in the set $F \cap B_\Psi(\rho,f^*)$. Theorem \ref{thm:small-ball-method} implies that with high probability, the empirical risk minimizer $\tilde{h}$ in $F \cap B_\Psi(\rho,f^*)$  satisfies
$$
\|\tilde{h}-h^*\|_{L_2(\mu)} \leq r(\rho).
$$
Moreover, it follows from \cite{LM13} and \cite{shahar-global-local} that under mild structural assumptions on $F$, $r(\rho)$ is the best possible error rate of any learning procedure in $F \cap B_\Psi(\rho,f^*)$ -- i.e., the minimax rate in that class.
\end{Remark}
Let ${\cal A}$ be the event from Corollary \ref{cor:small-ball-with-O} and set
$$
\gamma_{\cal O}(\rho) = \sup_{w \in \cal A} {\cal O}(\rho).
$$
$\gamma_{\cal O}$ will be of little importance in what follows, because it may be upper bounded by $(\theta/8)r^2(\rho)$. However, it will be of the utmost importance in \cite{LM_reg_comp}, where complexity-based regularization is studied (see Section \ref{sec:concluding-remarks} for more details).

\section{The main result} 
\label{sub:a_general_sparsity_error_bound}

Let us turn to the second part of the regularized functional -- namely, $\lambda(\Psi(f)-\Psi(f^*))$. Let $E^*$ be the dual space to $(E,\Psi)$ and set $\Psi^*$ to be the dual norm. $B_{\Psi^*}$ and $S_{\Psi^*}$ denote the dual unit ball and unit sphere, respectively; i.e., $B_{\Psi^*}$ consists of all the linear functionals $z^*$ on $E$ for which $\sup_{\Psi(x)=1} |z^*(x)| \leq 1$.
\begin{Definition}
The functional $z^* \in S_{\Psi^*}$ is a norming functional for $z \in E$ if $z^*(z)=\Psi(z)$.
\end{Definition}

In the language of Convex Analysis, a functional is norming for $x$ if and only if it belongs to $(\partial \Psi)_x$, the subdifferential of $\Psi$ in $x$.

\vskip0.4cm

Let $\Gamma_{f^*}(\rho)$ be the collection of functionals that are norming for some $f \in B_{\Psi}(\rho/20,f^*)$. In particular, $\Gamma_{f^*}(\rho)$ contains all the norming functionals of $f^*$.

Set
\begin{equation*}
\Delta(\rho) = \inf_{h\in H} \sup_{z^* \in \Gamma_{f^*}(\rho)} z^*(h-f^*),
\end{equation*}
where the infimum is taken in the set
$$
H = F \cap S_\Psi(\rho,f^*) \cap D(r(\rho),f^*) = \{ h \in F: \Psi(h-f^*)=\rho \ {\rm and} \ \|h-f^*\|_{L_2(\mu)} \leq r(\rho)\}.
$$

Note that if $z^* \in \Gamma_{f^*}(\rho)$ and $h \in S_{\Psi}(\rho,f^*)$ then $|z^*(h-f^*)| \leq \Psi(h-f^*)=\rho$. Thus, a lower bound of the form $\Delta(\rho) \geq (1-\delta)\rho$ implies that $\Gamma_{f^*}(\rho)$ is a relatively large subset of the dual unit sphere: each point in $F \cap S_\Psi(\rho,f^*) \cap D(r(\rho),f^*)$ has an `almost norming' functional in $\Gamma_{f^*}(\rho)$.

Our main result is that if $\Gamma_{f^*}(\rho)$ is indeed large enough to ensure that $\Delta(\rho) \geq 4/5 \rho$ then with high probability $\|\hat{f}-f^*\|_{L_2(\mu)} \leq r(\rho)$ and $\Psi(\hat f-f^*)\leq \rho$.

\begin{Theorem} \label{lemma:basic-combining-loss-and-reg}
Assume that $F$ is closed and convex. Let $\rho>0$ and set ${\cal A}$ to be an event on which Corollary \ref{cor:small-ball-with-O} holds. If $\Delta(\rho) \geq 4\rho/5$ and
$$
3\frac{\gamma_{\cal O}(\rho)}{\rho} \leq \lambda < \frac{\theta}{2} \cdot \frac{r^2(\rho)}{\rho},
$$
then on the event ${\cal A}$, a regularized empirical minimizer $\hat f \in {\rm argmin}_{f \in F} P_N {\cal L}_f^\lambda$ satisfies
$$
\Psi(\hat f-f^*)\leq \rho \mbox{ and } \|\hat{f}-f^*\|_{L_2(\mu)} \leq r(\rho).
$$
Moreover, since $r_{\cal O}(\rho) \leq (\theta/8)r^2(\rho)$, the same assertion holds if
\begin{equation*}
\frac{3\theta}{8} \cdot \frac{r^2(\rho)}{\rho} \leq \lambda < \frac{\theta}{2} \cdot \frac{r^2(\rho)}{\rho}.
\end{equation*}
\end{Theorem}

The proof of the theorem follows in three steps: first, one has to show that $P_N{\cal L}_f^\lambda$ is positive on the set $F \cap S_{\Psi}(\rho,f^*)$. Second, thanks to certain homogeneity properties of the functional, it is positive in $F \backslash B_{\Psi}(\rho,f^*)$, because it is positive on the `sphere' $F \cap S_\Psi(\rho,f^*)$. Finally, one has to study the functional in $F \cap B_{\Psi}(\rho,f^*)$ and verify that it is positive in that set, provided that $\|f-f^*\|_{L_2(\mu)} \geq r(\rho)$.

\begin{figure}[!h]
\centering
\begin{tikzpicture}[scale=0.2]
\draw (0,0) circle (8cm);
\draw (-2,0.3) node {$f^*$};
\filldraw (0,0) circle (0.2cm);
\draw (-10,0) -- (0,10) -- (10,0) -- (0,-10) -- (-10,0);
\draw (14.5,0) node {$S_\Psi(\rho, f^*)$};
\draw (12,6) node {$D(r(\rho),f^*)$};
\draw (0,0) -- (14,-8);
\draw (15,-8) node {$f$};
\filldraw (14,-8) circle (0.2cm);
\draw (4.7,-3.8) node {$h$};
\filldraw (6.42,-3.63) circle (0.2cm);
\end{tikzpicture}
\label{fig:partition_set_F}
\end{figure}

\vskip0.4cm

\proof Fix $h \in F \cap S_\Psi(\rho,f^*)$ and we shall treat two different cases: when $\|h-f^*\|_{L_2(\mu)} \geq r(\rho)$ and when $\|h-f^*\|_{L_2(\mu)} \leq r(\rho)$.

If $\|h-f^*\|_{L_2} \geq r(\rho)$, then by the triangle inequality for $\Psi$,
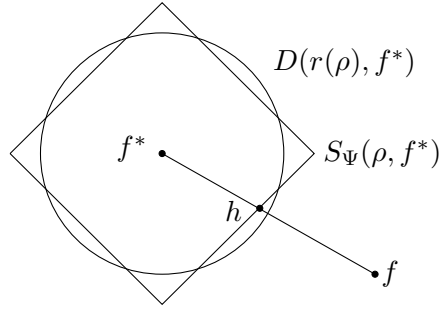
$$
\Psi(h)-\Psi(f^*)=\Psi(h-f^*+f^*)-\Psi(f^*) \geq -\Psi(h-f^*).
$$
Hence, for $(X_i,Y_i)_{i=1}^N \in {\cal A}$ and by the upper estimate in the choice of $\lambda$,
\begin{equation} \label{eq:large-distances-in-proof}
P_N {\cal L}_h^\lambda \geq \frac{\theta}{2} \|h-f^*\|_{L_2(\mu)}^2 - \lambda \Psi(h-f^*) \geq \frac{\theta}{2}r^2(\rho) - \lambda \rho > 0.
\end{equation}
Next, if $\|h-f^*\|_{L_2(\mu)} \leq r(\rho)$ then
\begin{equation*}
P_N {\cal L}_h^\lambda \geq -2{\cal O}(\rho) + \lambda(\Psi(h)-\Psi(f^*)).
\end{equation*}
Consider $u,v \in E$ that satisfy $f^*=u+v$ and $\Psi(u) \leq \rho/20$.
Let $z^*$ be any norming functional of $v$; thus, $z^* \in S_{\Psi^*}$ and $z^*(v)=\Psi(v)$. Since $\Psi(h)=\sup_{x^* \in B_{\Psi^*}} x^*(h)$ it follows that
\begin{align*}
\Psi(h)-\Psi(f^*) \geq & \Psi(h) - \Psi(v) -\Psi(u) \geq z^*(h-v) - \Psi(u)\geq  z^*(h-f^*) - 2 \Psi(u).
\end{align*}
This holds for any $v\in B_\Psi(\rho/20, f^*)$, and by the definition of $\Delta(\rho)$ and for an optimal choice of $z^*$,
\begin{equation} \label{eq:O-control-1}
P_N {\cal L}_h^\lambda \geq  -2{\cal O}(\rho) + \lambda(z^*(h-f^*) - 2 \Psi(u))
\geq  -2{\cal O}(\rho) + \lambda (\Delta(\rho) -\rho/10) >0,
\end{equation}
where the last inequality holds because $\Delta(\rho)\geq 4 \rho/5$ and $\lambda\geq 3\gamma_{\cal O}(\rho)/\rho$. Also, since $\gamma_{\cal O}(\rho) \leq (\theta/8)r^2(\rho)$, it suffices that $\lambda \geq (3\theta/8)r^2(\rho)/\rho$ to ensure that $P_N {\cal L}_h^\lambda >0$ in \eqref{eq:O-control-1}. This completes the proof of the first step -- that $P_N {\cal L}_h^\lambda >0$ on $F \cap S_{\Psi}(\rho,f^*)$.

\vskip0.4cm
Turning to the second step, one has to establish a similar inequality for functions outside $B_{\Psi}(\rho,f^*)$. To that end, let $f \in F \backslash B_\Psi(\rho,f^*)$. Since $F$ is convex and $\Psi$ is homogeneous, $f=f^*+\alpha (h-f^*)$ for some $h \in F \cap S_\Psi(\rho,f^*)$ and $\alpha>1$. Therefore,
$$
P_N {\cal Q}_{f-f^*} = \alpha^2 P_N {\cal Q}_{h-f^*}  {\mbox{ and }}  P_N {\cal M}_{f-f^*} = \alpha P_N {\cal M}_{h-f^*};
$$
moreover, $\Psi(f-f^*)=\alpha \Psi(h-f^*)$ and for every functional $z^*$,
$z^*(f-f^*)=\alpha z^*(h-f^*)$.

Thus, by \eqref{eq:large-distances-in-proof}, when $\|h-f^*\|_{L_2(\mu)} \geq r(\rho)$, $P_N {\cal L}_f^\lambda >0$, and when $\|h-f^*\|_{L_2(\mu)} \leq r(\rho)$,
\begin{align*}
P_N {\cal L}_f^\lambda & \geq \alpha^2 P_N {\cal Q}_{h-f^*} + 2\alpha P_N {\cal M}_{h-f^*} +  \lambda(\alpha z^*(h-f^*) - 2\Psi(u))\\
& \geq \alpha \bigl(P_N {\cal Q}_{h-f^*} + 2 P_N {\cal M}_{h-f^*} +  \lambda( z^*(h-f^*) - 2\Psi(u)) \bigr)>0.
\end{align*}

Finally, when $h \in F \cap B_{\Psi}(\rho,f^*)$ and $\norm{h-f^*}_{L_2(\mu)} \geq r(\rho)$, \eqref{eq:large-distances-in-proof} shows that $P_N{\cal L}_f^\lambda>0$.
\endproof

\begin{Remark} \label{rem:large-rho}
Note that if $\rho \geq \Psi(f^*)$ there is no upper limitation on the choice of $\lambda$. Indeed, if $\|h-f^*\|_{L_2(\mu)} \geq r(\rho)$ and $\Psi(h)=\rho \geq \Psi(f^*)$ then $\lambda(\Psi(h)-\Psi(f^*)) \geq 0$, and $P_N{\cal L}_h^\lambda >0$ just as in \eqref{eq:large-distances-in-proof}. The rest of the proof remains unchanged.
\end{Remark}

It follows from the proof that the quadratic component $P_N \cQ_{f-f^*}$ and the regularization one $\lambda(\Psi(f)-\Psi(f^*))$ dominate the multiplier component $2P_N\cM_{f-f^*}$ in different parts of $F$. The behaviour of $P_N {\cal Q}_{f-f^*}$ allows one to exclude the set $(F \cap B_{\psi}(\rho,f^*)) \backslash D(r(\rho),f^*)$, as well as any point in $F$ for which the interval $[f,f^*]$ intersects $(F \cap S_{\psi}(\rho,f^*)) \backslash D(r(\rho),f^*)$. This exclusion is rather free-of-charge, as it holds with no assumptions on the norm $\Psi$.

The situation is more subtle when trying to exclude points for which the interval $[f,f^*]$ intersects $F \cap S_{\psi}(\rho,f^*) \cap D(r(\rho),f^*)$.
That is precisely the region in which the specific choice of $\Psi$ is important and the regularization component is the reason why $P_N {\cal L}_f^\lambda>0$.

Figure~\ref{fig:partition_set_F_positive_excess_loss} shows this idea: $P_N {\cal L}_f^\lambda>0$ for two different reasons:  either $Q>M$ -- the quadratic component dominates the multiplier component, or $R>M$ -- the regularization component dominates the multiplier component.

Note  that an output of the sparsity equation is that the descent cone  $T_\Psi(f^*) = \cup_{\tau>0}\{h:\Psi(f^*+\tau h)\leq \Psi(f^*)\}$ does not intersect $S_\Psi(\rho, f^*)\cap D(r(\rho), f^*)$ when the ``sparsity condition'' $\Delta(\rho)\geq 4\rho/5$ is satisfied (cf. Figure~\ref{fig:descent_cone}).

\begin{figure}[h]
\centering
\begin{minipage}{.5\textwidth}
  \centering
\begin{tikzpicture}[scale=0.27]
\draw (0,0) circle (8cm);
\draw (2,0) node {$f^*$};
\filldraw (0,0) circle (0.2cm);
\draw (-10,0) -- (0,10) -- (10,0) -- (0,-10) -- (-10,0);
\draw[style = dashed] (-13,-4.2) -- (13,4.2);
\draw[style = dashed] (-13,4.2) -- (13,-4.2);
\draw[style = dashed] (-4.2,-13) -- (4.2,13);
\draw[style = dashed] (4.2,-13) -- (-4.2,13);
\filldraw[color=red, very thick, fill = red, opacity = 0.2] (4.2,13) -- (2.47,7.64) -- (7.64,2.47) -- (13,4.2);
\filldraw[color=red, very thick, fill = red, opacity = 0.2] (4.2,13) -- (13,13) -- (13,4.2);
\filldraw[color=red, very thick, fill = red, opacity = 0.2] (-4.2,-13) -- (-2.47,-7.64) -- (-7.64,-2.47) -- (-13,-4.2);
\filldraw[color=red, very thick, fill = red, opacity = 0.2] (-4.2,-13) -- (-13,-13) -- (-13,-4.2);
\filldraw[color=red, very thick, fill = red, opacity = 0.2] (-4.2,13) -- (-2.47,7.64) -- (-7.64,2.47) -- (-13,4.2);
\filldraw[color=red, very thick, fill = red, opacity = 0.2] (-4.2,13) -- (-13,13) -- (-13,4.2);
\filldraw[color=red, very thick, fill = red, opacity = 0.2] (4.2,-13) -- (2.47,-7.64) -- (7.64,-2.47) -- (13,-4.2);
\filldraw[color=red, very thick, fill = red, opacity = 0.2] (4.2,-13) -- (13,-13) -- (13,-4.2);
\filldraw[color=blue, very thick, fill = blue, opacity = 0.1] (-4.2,13) -- (-2.47,7.64) arc (110:70:7.2) -- (4.2,13) -- cycle;
\filldraw[color=blue, very thick, fill = blue, opacity = 0.1] (-4.2,-13) -- (-2.47,-7.64) arc (-110:-70:7.2) -- (4.2,-13);
\filldraw[color=blue, very thick, fill = blue, opacity = 0.1] (-13,4.2) -- (-7.64,2.47) arc (160:200:7.2) -- (-13,-4.2);
\filldraw[color=blue, very thick, fill = blue, opacity = 0.1] (13,4.2) -- (7.64,2.47) arc (20:-20:7.2) -- (13,-4.2);
\draw (8,8) node {$R>M$};
\draw (-8,8) node {$R>M$};
\draw (8,-8) node {$R>M$};
\draw (-8,-8) node {$R>M$};
\draw (0,11) node {$Q>M$};
\draw (-11,2) node {$Q>M$};
\draw (11,2) node {$Q>M$};
\draw (0,-11) node {$Q>M$};
\end{tikzpicture}
\captionof{figure}{The ``$Q>M$ and $R>M$'' decomposition.}
\label{fig:partition_set_F_positive_excess_loss}
\end{minipage}%
\begin{minipage}{.6\textwidth}
  \centering
\begin{tikzpicture}[scale=0.33]
      \tikzstyle{information text}=[rounded corners,fill=red!10,inner sep=1ex]
\filldraw [white!70!black] (0,0) circle (132pt);
\draw (7,4) node {$D(r(\rho), f^*)$};
\draw (2,0) node {$f^*$};
\filldraw (0,0) circle (0.2cm);
\filldraw[green] (-10,10) .. controls (-8,2) and (-8,2) .. (0,0) .. controls (-8,-2) and (-8,-2) .. (-10,-10) -- cycle;
\draw (-8,9) node {$T_\Psi(f^*)$};
\draw (0,10) .. controls (2,2) and (2,2) .. (10,0);
\draw (10,0) .. controls (2,-2) and (2,-2) .. (0,-10);
\draw (0,-10) .. controls (-2,-2) and (-2,-2) .. (-10,0);
\draw (-10,0) .. controls (-2,2) and (-2,2) .. (0,10);
\draw (3,9) node {$S_\Psi(\rho, f^*)$};
\draw [line width=3pt, blue] (1.6,4.4) .. controls (2.6,2.6) and (2.6,2.6) .. (4.4,1.6);
\draw [line width=3pt, blue] (1.6,-4.4) .. controls (2.6,-2.6) and (2.6,-2.6) .. (4.4,-1.6);
\draw [line width=3pt, blue] (-1.6,-4.4) .. controls (-2.6,-2.6) and (-2.6,-2.6) .. (-4.4,-1.6);
\draw [line width=3pt, blue] (-1.6,4.4) .. controls (-2.6,2.6) and (-2.6,2.6) .. (-4.4,1.6);
\end{tikzpicture}
\captionof{figure}{$T_\Psi(f^*)\cap S_\Psi(\rho, f^*)\cap D(r(\rho), f^*)=\emptyset$.}
\label{fig:descent_cone}
\end{minipage}
\end{figure}



\section{The role of $\Delta(\rho)$}
It is clear that $\Delta(\rho)$ plays a crucial role in the proof of Theorem \ref{lemma:basic-combining-loss-and-reg}, and that the larger $\Gamma_{f^*}(\rho)$ is, the better the lower bound on $\Delta(\rho)$.

Having many norming functionals of points in $B_\Psi(\rho/20,f^*)$ can be achieved somewhat artificially, by taking $\rho \sim \Psi(f^*)$. If $\rho$ is large enough, $B_\Psi(\rho/20,f^*)$ contains a $\Psi$-ball centred in $0$. Therefore, $\Gamma_{f^*}(\rho)$ is the entire dual sphere and $\Delta(\rho)=\rho$. This is the situation when one attempts to derive complexity-based bounds (see Section \ref{sec:concluding-remarks} and \cite{LM_reg_comp}),  i.e., when one wishes to find $\hat{f}$ that inherits some of $f^*$'s `good qualities' that are captured by $\Psi(f^*)$.

Here, we are interested in cases in which $\rho$ may be significantly smaller than $\Psi(f^*)$ and enough norming functionals have to be generated by other means.

If $\Psi$ is smooth, each $f \not = 0$ has a unique norming functional, and for a small $\rho$, the norming functionals of points in $B_\Psi(\rho/20,f^*)$ are close to the (unique) norming functional of $f^*$; hence there is little hope that $\Gamma_{f^*}(\rho)$ will be large enough to ensure that $\Delta(\rho) \sim \rho$. It is therefore reasonable to choose $\Psi$ that is not smooth in $f^*$ or in a neighbourhood of $f^*$.

Another important fact is that $\Gamma_{f^*}(\rho)$ need not be as large as the entire dual sphere to ensure that $\Delta(\rho) \sim \rho$. Indeed, it suffices if $\Gamma_{f^*}(\rho)$ contains `almost norming' functionals only to points that satisfy $\|w\|_{L_2(\mu)} \leq r(\rho)/\rho$ and $\Psi(w)=1$, rather than to every point in the sphere $S_\Psi$.

\subsection{$\Delta(\rho)$ and sparsity}
It turns out that the combination of the right notion of sparsity with a wise choice of a norm $\Psi$ ensures that $\Gamma_{f^*}(\rho)$ contains enough `almost norming' functionals precisely for the subset of the sphere one is interested in.

To give an indication of how this happens, let us show the following:

\begin{Lemma} \label{lemma:support-property}
Let $Z \subset S_{\Psi^*}$, $W \subset S_\Psi$ and $0<\eta_1,\eta_2<1$. If every $w \in W$ can be written as $w=w_1+w_2$, where $\Psi(w_1) \leq \eta_1 \Psi(w)$ and $\sup_{z^* \in Z} z^*(w_2) \geq (1-\eta_2) \Psi(w_2)$, then
$$
\inf_{w \in W} \sup_{z^* \in Z} z^*(w) \geq (1-\eta_1)(1-\eta_2)-\eta_1
$$
In particular, if $\eta_1,\eta_2 \leq 1/20$ then $\inf_{w \in W} \sup_{z^* \in Z} z^*(w) \geq 4/5$.
\end{Lemma}

\proof
Let $w=w_1+w_2$ and observe that $\Psi(w_2) \geq \Psi(w)-\Psi(w_1) \geq (1-\eta_1)\Psi(w)$.
Thus, for the optimal choice of $z^* \in Z$,
\begin{align*}
z^*(w_1+w_2) \geq &(1-\eta_2)\Psi(w_2) + z^*(w_1) \geq (1-\eta_2)\Psi(w_2) -\eta_1 \Psi(w).
\\
\geq  & \bigl((1-\eta_1)(1-\eta_2)-\eta_1\bigr) \Psi(w),
\end{align*}
and the claim follows because $w \in S_{\Psi}$.
\endproof

Let $E=\R^d$ viewed as a class of linear functionals on $\R^d$. Set $\mu$ to be an isotropic measure on $\R^d$; thus $\{t \in \R^d : \E \inr{t,X}^2 \leq 1\}=B_2^d$.

Assume that for $t \in \R^d$ that is supported on $I \subset \{1,...,d\}$, the set of its norming functionals consists of functionals of the form $z_0^*+(1-\eta_2)u^*$ for some fixed $z_0^*$ that is supported on $I$ and {\it any} $u \in B_{\Psi^*}$ that is supported on $I^c$ (such is the case, for example, when $E=\ell_1^d$).

For every such $t$, consider $w \in \rho S_{\Psi}$ and set $w_1=P_I w$ and $w_2 = P_{I^c} w$, the coordinate projections of $w$ onto ${\rm span}(e_i)_{i \in I}$ and ${\rm span}(e_i)_{i \in I^c}$, respectively. Hence, there is a functional $z^*=z_0^*+(1-\eta_2)u^*$ that is norming for $t$ and also satisfies
$$
z^*(w_2) = (1-\eta_2)u^*(w_2)=(1-\eta_2)\Psi(w_2).
$$
Therefore, Lemma \ref{lemma:support-property} may be applied once $\Psi(P_I w) \leq \eta_1 \Psi(w)$.

Naturally, such a shrinking phenomenon need not be true for {\it every} $w \in S_{\Psi}$; fortunately, it is only required for $w \in S_\Psi \cap (r(\rho)/\rho)D$ -- and we will show that it is indeed the case in the three examples we present. In all three, the combination of sparsity and the right choice of the norm helps in establishing a lower bound on $\Delta(\rho)$ in two ways: firstly, the set $\Gamma_{t^*}(\rho)$ consists of functionals that are `almost norming' for any $x$ whose support is disjoint from the support of $t^*$; and secondly, a coordinate projection `shrinks' the $\Psi$ norm of points in $\rho S_\Psi \cap r(\rho)D$.

\subsection{$\Delta(\rho)$ in the three examples}
Let us show that in the three examples, the LASSO, SLOPE and trace norm regularization, $\Delta(\rho) \geq (4/5)\rho$ for the right choice of $\rho$, and that choice depends on the degree of sparsity in each case.

In all three examples, we will assume that the underlying measure is isotropic; thus the $L_2(\mu)$ norm coincides with the natural Euclidean structure: the $\ell_2^d$ norm for the LASSO and SLOPE, and the Hilbert-Schmidt norm for trace-norm regularization.

\vskip0.3cm

\noindent{\bf The LASSO.}

Observe that if $f^*=\inr{t^*,\cdot}$ is the true minimizer of the functional $\inr{t,\cdot}\to \E(\inr{t,X}-Y)^2$ in $F=\{\inr{t, \cdot}:t\in\R^d\}$, then any function $h_t=\inr{t,\cdot}$ for which $\|h_t-f^*\|_{L_2} \leq r(\rho)$ and $\Psi(h_t-f^*) = \rho$ is of the form $h_t=\inr{t,\cdot}=\inr{w+t^*,\cdot}$, where $w \in \rho S(\ell_1^d) \cap r(\rho)B_2^d.$ Recall that the dual norm to $\|\cdot\|_1$ is $\|\cdot\|_\infty$, and thus
$$
\Delta(\rho) = \inf_{w \in \rho S(\ell_1^d) \cap r(\rho)B_2^d} \sup_{z \in \Gamma_{t^*}(\rho)}  \inr{z,w},
$$
where $\Gamma_{t^*}(\rho)$ is the set of all vectors $z^*\in \R^d$ that satisfy
$$
\|z^*\|_\infty = 1 \ \ {\rm and} \ \ z^*(v)=\|v\|_1 \ {\rm for \ some \ } v \ {\rm for \ which \ } \|v-t^*\|_1 \leq \rho/20.
$$

\begin{Lemma} \label{lemma:norming-for-LASSO}
If $t^*=v+u$ for $u \in (\rho/20)B_1^d$ and $100|{\rm supp}(v)| \leq (\rho/r(\rho))^2$ then $\Delta(\rho) \geq 4\rho/5$.
\end{Lemma}
In other words, if $t^*$ is well approximated with respect to the $\ell_1^d$ norm by some $v \in \R^d$ that is $s$-sparse, and $s$ is small enough relative to the ratio $(\rho/r(\rho))^2$, then $\Delta(\rho) \geq (4/5)\rho$.

\vskip0.4cm

Just as noted earlier, we shall use two key properties of the $\ell_1$ norm and sparse vectors: firstly, that if $x$ and $y$ have disjoint supports, there is a functional that is simultaneously norming for $x$ and $y$, i.e., $z^* \in B_\infty^d$ for which
\begin{equation} \label{eq:sim-norm-ell-1}
z^*(x)=\|x\|_1 \ \ {\rm and} \ \ z^*(y)=\|y\|_1;
\end{equation}
secondly, that if $\|x\|_1 = \rho$ and $\|x\|_2$ is significantly smaller than $\rho$, a coordinate projection `shrinks' the $\ell_1^d$ norm: $\|P_I x\|_1$ is much smaller than $\|x\|_1$.

\proof
Let $w \in \rho S(\ell_1^d) \cap r(\rho)B_2^d$. Since $\norm{t^*-v}_1\leq \rho/20$ there exists $z^*\in \Gamma_{t^*}(\rho)$ that is norming for $v$. Moreover, if $I={\rm supp}(v)$, then according to \eqref{eq:sim-norm-ell-1} one can choose $z^*$ that is also norming for $P_{I^c}w$. Thus, $\norm{P_{I^c}w}_1 = z^*(P_{I^c}w)$ and
\begin{align*}
z^*(w) = z^*(P_I w) + z^*(P_{I^c}w)\geq \norm{P_{I^c} w}_1 - \norm{P_I w}_1\geq \norm{w}_1 - 2 \norm{P_I w}_1.
\end{align*}
Since $\norm{w}_2\leq r(\rho)$, one has
$\norm{P_I w}_1\leq \sqrt{s} \norm{P_I w}_2\leq \sqrt{s}r(\rho)$. Therefore,
$$
\inr{z,w}\geq \rho - 2\sqrt{s}r(\rho)\geq 4\rho/5
$$
when $100 s \leq (\rho/r(\rho))^2$.
\endproof

\vskip0.4cm

\noindent{\bf SLOPE.}

Let $\beta_1 \geq \beta_2 \geq ... \geq \beta_d>0$ and recall that $\Psi(t)=\sum_{i=1}^d \beta_i t_i^*$.

Note that $\Psi(t) = \sup_{z \in Z} \inr{z,t}$, for
$$
Z=\left\{ \sum_{i=1}^d \eps_i \beta_{\pi_i} e_i : \ (\eps_i)_{i=1}^d \in \{-1,1\}^d, \ \pi \ {\rm is \ a \ permulation \ of \ } \{1,...,d\} \right\}.
$$
Therefore, the extreme points of the dual unit ball are of the form $\sum_{i=1}^d \eps_i \beta_{\pi_i} e_i$.

Following the argument outlined above, let us show that if $x$ is supported on a reasonably small $I \subset \{1,...,d\}$, the set of norming functionals of $x$ consists of `almost norming' functionals for any $y$ that is supported on $I^c$. Moreover, and just like the $\ell_1^d$ norm, if $\Psi(x)=\rho$ and $\|x\|_2$ is significantly smaller than $\rho$, a coordinate projection of $x$ `shrinks' its $\Psi$ norm.

\begin{Lemma} \label{lemma:delta-for-slope}
Let $1 \leq s \leq d$ and set ${\cal B}_s = \sum_{i \leq s} \beta_i/\sqrt{i}$. If $t^*$ is $\rho/20$ approximated (relative to $\Psi$) by an $s$-sparse vector and if $40{\cal B}_s \leq \rho/r(\rho)$ then $\Delta(\rho) \geq 4\rho/5$.
\end{Lemma}

\proof Let $t^*=u+v$, for $v$ that is supported on at most $s$ coordinates and $u \in (\rho/20)B_\Psi$. Set $I \subset \{1,...,d\}$ to be the support of $v$ and let $z=(z_i)_{i=1}^d$ be a norming functional for $v$ to be specified later; thus, $z \in \Gamma_{t^*}(\rho)$.

Given $t$ for which $\Psi(t-t^*)=\rho$ and $\|t-t^*\|_2 \leq r(\rho)$, one has
\begin{align*}
z(t-t^*) & = z(t-v) - z(u) = z(P_{I^c}(t-v)) + z(P_I(t-v)) - z(u)
\\
\geq & \sum_{i \in I^c} z_i (t-v)_i + \sum_{i \in I} z_i (t-v)_i - \Psi(u)
\\
\geq & \sum_{i \in I^c} z_i (t-v)_i - \sum_{i \leq s} \beta_i (t-v-u)_i^\sharp - 2\Psi(u)
\\
 = & \sum_{i \in I^c} z_i (t-v)_i - \sum_{i \leq s} \beta_i (t-t^*)_i^\sharp - 2\Psi(u) = (*).
\end{align*}
Since $v$ is supported in $I$, one may optimize the choice of $z$ by selecting the right permutation of the coordinates in $I^c$, and
\begin{align*}
& \sum_{i \in I^c} z_i (t-v)_i \geq \sum_{i >s} \beta_i (t-v)_i^\sharp \geq \sum_{i > s} \beta_i(t-v-u)_i^\sharp - \Psi(u)
\\
 = & \sum_{i=1}^d \beta_i(t-t^*)_i^\sharp - \sum_{i \leq s} \beta_i (t-t^*)_i^\sharp - \Psi(u).
\end{align*}
Therefore,
\begin{equation*}
(*) \geq \sum_{i=1}^d \beta_i(t-t^*)_i^\sharp - 2\sum_{i \leq s} \beta_i (t-t^*)_i^\sharp - 3\Psi(u) \geq \frac{17}{20}\rho - 2\sum_{i \leq s} \beta_i (t-t^*)_i^\sharp.
\end{equation*}

Since $\|t-t^*\|_2 \leq r(\rho)$, it is evident that $(t-t^*)_i^\sharp \leq r(\rho)/\sqrt{i}$, and
$$
\sum_{i=1}^s \beta_i (t-t^*)_i^\sharp \leq r(\rho) \sum_{i=1}^s \frac{\beta_i}{\sqrt{i}} = r(\rho) {\cal B}_s.
$$
Hence, if $\rho \geq 40 r(\rho) {\cal B}_s$ then $\Delta(\rho) \geq 4\rho/5$.
\endproof

\vskip0.4cm
\noindent{\bf Trace-norm regularization.}

The trace norm has similar properties to the $\ell_1$ norm. Firstly, one may show that the dual norm to $\| \cdot \|_1$ is $\| \cdot \|_\infty$, which is simply the standard operator norm. Moreover, one may find a functional that is simultaneously norming for any two elements with `disjoint support' (and of course, the meaning of `disjoint support' has to be interpreted correctly here).
Finally, it satisfies a `shrinking' phenomenon for matrices whose Hilbert-Schmidt norm is significantly smaller than their trace norm.
\begin{Lemma} \label{lemma:norming-for-trace-norm}
If $A^*=V+U$, where  $\|U\|_1 \leq \rho/20$ and $400{\rm rank}(V) \leq (\rho/r(\rho))^2$, then $\Delta(\rho) \geq 4\rho/5$.
\end{Lemma}

The fact that a low-rank matrix has many norming functionals is  well known and follows, for example, from \cite{MR1160950}.

\begin{Lemma} \label{lemma:low-rank-matrix}
Let $V \in R^{m \times T}$ and assume that $V=P_I V P_J$ for appropriate orthogonal projections onto subspaces $I \subset \R^m$ and $J \subset \R^T$. Then, for every $W \in \R^{m \times T}$ there is a matrix $Z$ that satisfies $\|Z\|_\infty=1$, and
$$
\inr{Z,V} = \|V\|_1,  \ \ \inr{Z,P_{I^\perp} W P_{J^\perp}}= \|P_{I^\perp} W P_{J^\perp}\|_1,
$$
$$
\inr{Z,P_{I} W P_{J^\perp}}=0 \ \ {\rm and} \ \ \inr{Z,P_{I^\perp} W P_{J}}=0.
$$
\end{Lemma}
Lemma \ref{lemma:low-rank-matrix} describes a similar phenomenon to the situation in $\ell_1^d$, but with a different notion of `disjoint support': if $V$ is low-rank and the projections $P_I$ and $P_J$ are non-trivial, one may find a functional that is norming both for $V$ and for the part of $W$ that is `disjoint' of $V$. Moreover, the functional vanishes on the `mixed' parts $P_{I} W P_{J^\perp}$ and $P_{I^\perp} W P_{J}$.

\vskip0.4cm

\noindent{\bf Proof of Lemma \ref{lemma:norming-for-trace-norm}.}
Recall that $S_1$ is the unit sphere of the trace norm and that $B_2$ is the unit ball of the Hilbert-Schmidt norm. Hence,
\begin{equation*}
\Delta(\rho) = \inf_{W\in \rho S_1 \cap r(\rho)B_2} \sup_{Z\in \Gamma_{A^*}(\rho)}\inr{Z,W}
\end{equation*}where $\Gamma_{A^*}(\rho)$ is the set of all matrices $Z\in \R^{m\times T}$ that satisfy $\norm{Z}_\infty = 1$ and $\inr{Z,V}=1$ for some $V$ for which $\norm{A^*-V}_1\leq \rho/20$.

Fix a rank-$s$ matrix $V=P_I V P_J$, for orthogonal projections $P_I$ and $P_J$ that are onto subspaces of dimension $s$. Consider $W\in \R^{m\times T}$ for which $\|W\|_1 = \rho$ and $\|W\|_2 \leq r(\rho)$ and put $Z$ to be a norming functional of $V$ as in Lemma \ref{lemma:low-rank-matrix}. Thus, $Z \in \Gamma_{A^*}(\rho)$ and
\begin{align*}
\inr{Z,W} = & \inr{Z,P_{I^\perp} W P_{J^\perp}} + \inr{Z,P_{I} W P_{J}} =\|P_{I^\perp} W P_{J^\perp}\|_1 - \|P_{I} W P_{J}\|_1  \nonumber
\\
\geq & \|W\|_1 - \|P_{I} W P_{J^\perp}\|_1 -\|P_{I^\perp} W P_{J}\|_1-2\|P_{I} W P_{J}\|_1.
\end{align*}
All that remains is to estimate the trace norms of the three components that are believed to be `low-dimension' - in the sense that their rank is at most $s$.

Recall that $(\sigma_i(A))$ are the singular values of $A$ arranged in a non-increasing order. It is straightforward to verify (e.g., using the characterization of the singular values via low-dimensional approximation), that
$$
\sigma_i(P_{I} W P_{J^\perp}),\sigma_i(P_{I^\perp} W P_{J}),\sigma_i(P_{I} W P_{J}) \leq \sigma_i(W).
$$
Moreover,  $\|W\|_2 \leq r(\rho)$, therefore, being rank-$s$ operators, one has
\begin{align*}
&\|P_{I} W P_{J^\perp}\|_1, \ \|P_{I^\perp} W P_{J}\|_1, \ \|P_{I} W P_{J}\|_1 \leq \sum_{i =1}^s \sigma_i(W)\leq \sqrt{s} \Big(\sum_{i=1}^s \sigma_i^2(W) \Big)^{1/2}\leq \sqrt{s}r(\rho),
\end{align*}
implying that
$$
\inr{Z,W} \geq \rho - 4r(\rho) \sqrt{s}.
$$
Therefore, if $400 s \leq (\rho/r(\rho))^2$, then $\Delta(\rho) \geq 4\rho/5$.
\endproof

\section{The three examples revisited} \label{sec:r}
The estimates on $\Delta(\rho)$ presented above show that in all three examples, when $f^*$ is well approximated by a function whose `degree of sparsity' is $\lesssim (\rho/r(\rho))^2$, then $\Delta(\rho) \geq 4\rho/5$ and Theorem \ref{lemma:basic-combining-loss-and-reg} may be used. Clearly, the resulting error rates depend on the right choice of $\rho$, and thus on $r(\rho)$.

Because $r(\rho)$ happens to be the minimax rate of the learning problem in the class $F\cap B_\Psi(\rho,f^*)$, its properties have been studied extensively. Obtaining an estimate on $r(\rho)$ involves some assumptions on $X$ and $\xi$, and the one setup in which it can be characterized for an arbitrary class $F$ is when the class is $L$-subgaussian and $\xi \in L_q$ for some $q>2$ (though $\xi$ need not be independent of $X$). It is straightforward to verify that an $L$-subgaussian class satisfies the small-ball condition of Assumption~\ref{ass:small-ball} for $\kappa=1/2$ and $\eps=c/L^4$ where $c$ is an absolute constant. Moreover, if the class is $L$-subgaussian, the natural complexity parameter associated with it is the expectation of the supremum of the canonical Gaussian process indexed by the class.

\begin{Definition}\label{def:gaussian_mean_width}
 Let $F \subset L_2(\mu)$ and set $\{G_f: f\in F\}$ to be the canonical Gaussian process indexed by $F$; that is, each $G_f$ is a centred Gaussian variable and the covariance structure of the process is endowed by the inner product in $L_2(\mu)$. The expectation of the supremum of the process is defined by
 \begin{equation*}
  \ell_*(F) = \sup\{ \E \sup_{f\in F^\prime} G_f : F^\prime \subset F \ {\rm is \ finite} \}.
  \end{equation*}
 \end{Definition}

It follows from a standard chaining argument that if $F$ is $L$-subgaussian then
\begin{equation*}
\E \sup_{f\in F}\Big|\frac{1}{N}\sum_{i=1}^N \eps_i (h-f^*)(X_i)\Big|\lesssim L\frac{ \ell_*(F)}{\sqrt{N}}.
\end{equation*}
Therefore, if
$$
F_{\rho,r} = F\cap B_\Psi(\rho,f^*)\cap D(r,f^*)
$$
then for every $\rho>0$ and $f^* \in F$
\begin{equation*}
r_Q(\rho) \leq \inf \left\{r>0 : \ell_*(F_{\rho,r})\leq C(L) r \sqrt{N} \right\}.
\end{equation*}
Turning to $r_M$, we shall require the following fact from \cite{shahar_multi_pro}.
\begin{Theorem}[Corollary~1.10 in \cite{shahar_multi_pro}]\label{thm:multi-main}
Let $q>2$ and $L \geq 1$. For every $0<\delta<1$ there is a constant $c=c(\delta,L,q)$ for which the following holds. If $H$ is an $L$-subgaussian class and $\xi \in L_q$, then with probability at least $1-\delta$,
\begin{equation*}
\sup_{h \in H} \left|\frac{1}{\sqrt{N}}\sum_{i=1}^N \eps_i\xi_i h(X_i)\right| \leq c \|\xi\|_{L_q}\ell_*(H).
\end{equation*}
\end{Theorem}
The complete version of Theorem \ref{thm:multi-main} includes a sharp estimate on the constant $c$. However, obtaining accurate probability estimates is not the main feature of this note and deriving such estimates leads to a cumbersome presentation. To keep our message to the point, we have chosen not to present the best possible probability estimates in what follows.

\vskip0.4cm

A straightforward application of Theorem \ref{thm:multi-main} shows that
\begin{equation*}
r_M(\rho) \leq  \inf \left\{r>0 : \norm{\xi}_{L_q}\ell_*(F_{\rho,r}) \leq c r^2 \sqrt{N} \right\}
\end{equation*}
for a constant $c$ that depends on $L,q$ and $\delta$.

\vskip0.4cm
Recall that we have assumed that $X$ is isotropic, which means that the $L_2(\mu)$ norm coincides with the natural Euclidean structure on the space: the standard $\ell_2^d$ norm for the LASSO and SLOPE and the Hilbert-Schmidt norm for trace norm regularization. Since the covariance structure of the indexing Gaussian process is endowed by the inner product, it follows that
$$
\ell_*(\rho B_\Psi \cap r D) = \E \sup_{w \in \rho B_\Psi \cap r B_2} \inr{G,w}
$$
for the standard Gaussian vector $G=(g_1,...,g_d)$ in the case of the LASSO and SLOPE and the Gaussian matrix $G=(g_{ij})$ in the case of trace norm minimization. Hence, one may obtain a bound on $r(\rho)$ by estimating this expectation in each case.

\vskip0.4cm
\noindent{\bf The LASSO and SLOPE.}
Let $(\beta_i)_{i=1}^d$ be a non-increasing positive sequence and set
$\Psi(t)=\sum_{i=1}^d t_i^\sharp \beta_i$.

Since the LASSO corresponds to the choice of $(\beta_i)_{i=1}^d=(1,...,1)$, it suffices to identify $\ell_*(\rho B_\Psi \cap r B_2^d)$ for the SLOPE norm and a general choice of weights.
\begin{Lemma} \label{lemma:monotone-for-slope}
There exists an absolute constant $C$ for which the following holds. If $\beta$ and $\Psi$ are as above, then
$$
\E \sup_{w \in \rho B_\Psi \cap r B_2^d} \inr{G,w} \leq C\min_k  \left\{r \sqrt{(k-1)\log\Big(\frac{ed}{k-1}\Big)} + \rho \max_{i \geq k} \frac{\sqrt{\log(ed/i)}}{\beta_i} \right\}
$$
(and if $k=1$, the first term is set to be $0$).
\end{Lemma}

\proof Fix $1 \leq k \leq d$. Let $J$ be the set of indices of the $k$ largest coordinates of $(|g_i|)_{i=1}^d$, and for every $w$ let $I_w$ be the sets of indices of the $k$ largest coordinates of $(|w_i|)_{i=1}^d$. Put $J_w=J \cup I_w$ and note that $|J_w| \leq 2k$. Hence,
\begin{align*}
 \sup_{w \in \rho B_\Psi \cap r B_2^d}& \sum_{i=1}^d w_i g_i \leq \sup_{w \in r B_2^d} \sum_{i \in J_w} w_i g_i + \sup_{w \in \rho B_{\Psi}} \sum_{i \in J_w^c} w_i g_i \nonumber
\\
&\lesssim  r \left(\sum_{i<k} (g_i^\sharp)^2 \right)^{1/2} + \sup_{w \in \rho B_{\Psi}} \sum_{i \geq k} w_i^\sharp \beta_i \frac{g_i^\sharp}{\beta_i} \lesssim r \left(\sum_{i<k} (g_i^\sharp)^2 \right)^{1/2} + \rho \max_{i \geq k} \frac{g_i^\sharp}{\beta_i}.
\end{align*}

As a starting point, note that a standard binomial estimate shows that
\begin{align*}
Pr\left(g_i^\sharp \geq t\sqrt{\log(ed/i)}\right) \leq & \binom{d}{i}Pr^i\left(|g| \geq t\sqrt{\log(ed/i)}\right) \nonumber
\\
\leq & 2\exp(i\log(ed/i)-i\log(ed/i) \cdot t^2/2).
\end{align*}
Applying the union bound one has that for $t \geq 4$, with probability at least $1-2\exp(-(t^2/2) k \log (ed/k))$,
\begin{equation}\label{eq:monotone-in-proof-1}
g_i^\sharp \leq c_3t \sqrt{\log(ed/i)} \ \ {\rm for \ every \ } i \geq k.
\end{equation}
The same argument shows that $\E (g_i^\sharp)^2 \lesssim \log(ed/i)$.

Let $U_k$ be the set of vectors on the Euclidean sphere that are supported on at most $k$ coordinates. Set
$$
\|x\|_{[k]} = \Bigl(\sum_{i \leq k} (x_i^\sharp)^2\Bigr)^{1/2}=\sup_{u \in U_k} \inr{x,u}
$$
and recall that by the Gaussian concentration of measure theorem (see, e.g., Theorem~7.1 in \cite{MR1849347}),
$$
\left(\E \|G\|_{[k]}^q\right)^{1/q} \leq \E \|G\|_{[k]} + c\sqrt{q} \sup_{u \in U_k} \|\inr{G,u}\|_{L_2} \leq \E\|G\|_{[k]}+c_1\sqrt{q}.
$$
Moreover, since $\E (g_i^\sharp)^2 \lesssim \log(ed/i)$, one has
$$
\E\|G\|_{[k]} \leq \Bigl(\E \sum_{i \leq k} (g_i^\sharp)^2\Bigr)^{1/2} \lesssim \sqrt{k \log(ed/k)}.
$$
Therefore, by Chebyshev's inequality for $q \sim k\log(ed/k)$, for $t \geq 1$, with probability at least $1-2t^{-c_1k\log(ed/k)}$,
$$
\Bigl(\sum_{i \leq k} (g_i^\sharp)^2\Bigr)^{1/2} \leq c_2 t \sqrt{k\log(ed/k)}.
$$
Turning to the `small coordinates', by \eqref{eq:monotone-in-proof-1},
$$
\max_{i \geq k} \frac{g_i^\sharp}{\beta_i} \lesssim t \max_{i \geq k} \frac{\sqrt{\log(ed/i)}}{\beta_i}.
$$
It follows that for every choice of $1 \leq k \leq d$,
\begin{align*}
\E \sup_{w \in \rho B_\Psi \cap r B_2^d} &\inr{G,w} \lesssim  r \E \Bigl(\sum_{i<k} (g_i^\sharp)^2 \Bigr)^{1/2} + \rho \E \max_{i \geq k} \frac{g_i^\sharp}{\beta_i} \nonumber
\\
\lesssim & r \sqrt{(k-1)\log(ed/(k-1))} + \rho \max_{i \geq k} \frac{\sqrt{\log(ed/i)}}{\beta_i},
\end{align*}
and, if $k=1$, the first term is set to be $0$.
\endproof

If $\beta=(1,...,1)$ (which corresponds to the LASSO), then $B_\Psi=B_1^d$, and one may select $\sqrt{k} \sim \rho/r$, provided that $r \leq \rho \leq r \sqrt{d}$. In that case,
\begin{equation*}
\E \sup_{w \in \rho B_1^d \cap r B_2^d} \inr{G,w} \lesssim \rho \sqrt{\log(edr^2/\rho^2)}.
\end{equation*}
The estimates when $r \geq \rho$ or $r\sqrt{d} \leq \rho$ are straightforward. Indeed, if $r \geq \rho$ then $\rho B_1^d \subset r B_2^d$ and
$$
\ell_*(\rho B_1^d\cap r B_2^d) = \ell_*(\rho B_1^d) \sim \rho \sqrt{\log(ed)},
$$
while if  $r\sqrt{d}\leq \rho$ then $r B_2^d \subset \rho B_1^d$, and
$$
\ell_*(\rho B_1^d\cap r B_2^d) = \ell_*(r B_2^d) \sim r \sqrt{d}.
$$
\vskip0.4cm

\noindent{\bf The LASSO}.

A straightforward computation shows that
\begin{equation*}
r_M^2(\rho) \lesssim_{L,q,\delta} \left\{
\begin{array}{cc}
\frac{\norm{\xi}_{L_q}^2 d}{N} & \mbox{ if } \rho^2 N \gtrsim_{L,q,\delta} \norm{\xi}_{L_q}^2 d^2
\\
\\
\rho \norm{\xi}_{L_q}\sqrt{\frac{1}{N}\log\Big(\frac{e\norm{\xi}_{L_q} d}{\rho\sqrt{N}}\Big)} & \mbox{ otherwise,}
\end{array}
\right.
\end{equation*}
and
\begin{equation*}
r_Q^2(\rho) \lesssim_L  \left\{
\begin{array}{cc}
  0 & \mbox{ if }  N \gtrsim_L  d \\
\frac{\rho^2}{N}\log\Big(\frac{c(L)d}{N}\Big) & \mbox{ otherwise}.
\end{array}
\right.
\end{equation*}

\noindent {\bf Proof of Theorem \ref{thm:intro-LASSO-est}.}
We will actually prove a slightly stronger result, which gives an improved estimation error if one has prior information on the degree of sparsity.

\vskip0.3cm

Using the estimates on $r_M$ and $r_Q$, it is straightforward to verify that the sparsity condition of Lemma \ref{lemma:norming-for-LASSO} holds when $N \gtrsim_{L,q,\delta} s \log(ed/s)$ and for any
 \begin{equation*}
  \rho \gtrsim_{L,q,\delta} \|\xi\|_{L_q} s \sqrt{\frac{1}{N}\log\Big(\frac{ed}{s}\Big)}.
  \end{equation*}
It follows from Lemma~\ref{lemma:norming-for-LASSO} that if there is an $s$-sparse vector that belongs to $t^*+(\rho/20)B_1^d$, then $\Delta(\rho)\geq 4\rho/5$. Finally, Theorem \ref{lemma:basic-combining-loss-and-reg} yields the stated bounds on $\|\hat{t}-t^*\|_1$ and $\|\hat{t}-t^*\|_2$ once we set
  \begin{equation*}
  \lambda \sim \frac{r^2(\rho)}{\rho} \sim_{L,q,\delta} \|\xi\|_{L_q} \sqrt{\frac{1}{N}\log\Big(\frac{ed}{s}\Big)}.
  \end{equation*}
The estimates on $\|\hat{t}-t^*\|_p$ for $1\leq p\leq 2$ can be easily verified because
$$
\norm{x}_p\leq \norm{x}_1^{-1+2/p}\norm{x}_2^{2-2/p}.
$$
In case one has no prior information on $s$, one may take
$$
\rho \sim_{L,q,\delta} \|\xi\|_{L_q} s \sqrt{\frac{1}{N}\log(ed)}
$$
and
\begin{equation*}
  \lambda \sim_{L,q,\delta} \|\xi\|_{L_q} \sqrt{\frac{\log(ed)}{N}}.
\end{equation*}
The rest of the argument remains unchanged.
\endproof

\vskip0.4cm

\noindent {\bf SLOPE}

Assume that $\beta_i \leq C\sqrt{\log(ed/i)}$, which is the standard assumption for SLOPE \cite{slope1,slope2}. By considering the cases $k=1$ and $k=d$,
\begin{equation}\label{eq:gauss_mean_slope_iso}
\E \sup_{w \in \rho B_\Psi \cap r B_2^d} \inr{G,w} \lesssim \min\{C\rho,\sqrt{d}r\}.
\end{equation}
Thus, one may show that
\begin{equation*}
r_Q^2(\rho) \lesssim_{L} \left\{
\begin{array}{cc}
0 & \mbox{ if }  N \gtrsim_L d
\\
\\
\frac{\rho^2}{N} & \mbox{ otherwise,}
\end{array}
\right.
\mbox{ and } \ \
r_M^2(\rho) \lesssim_{L,q,\delta} \left\{
\begin{array}{cc}
\|\xi\|_{L_q}^2 \frac{ d}{N} & \mbox{ if } \rho^2 N \gtrsim_{L,q,\delta} \|\xi\|_{L_q}^2 d^2
\\
\\
\|\xi\|_{L_q}\frac{\rho}{\sqrt{N}} & \mbox{ otherwise.}
\end{array}
\right.
\end{equation*}

\noindent{\bf Proof of Theorem \ref{thm:intro-SLOPE-est}.}
Recall that $\cB_s = \sum_{i\leq s}\beta_i/\sqrt{i}$, and when $\beta_i\leq C \sqrt{\log(ed/i)}$, one may verify that
$$
\cB_s\lesssim C \sqrt{s \log(ed/s)}.
$$
Hence, the condition $\cB_s\lesssim \rho/r(\rho)$ holds when $N \gtrsim_{L,q,\delta} s \log(ed/s)$ and
\begin{equation*}
\rho \gtrsim_{L,q,\delta} \|\xi\|_{L_q} \frac{s }{\sqrt{N}} \log\Big(\frac{ed}{s}\Big).
\end{equation*}
It follows from Lemma~\ref{lemma:delta-for-slope} that $\Delta(\rho)\geq 4\rho/5$ when there is an $s$-sparse vector in $t^*+(\rho/20) B_\Psi$; therefore, one may apply Theorem \ref{lemma:basic-combining-loss-and-reg} for the choice of
\begin{equation*}
\lambda \sim \frac{r^2(\rho)}{\rho} \sim_{L,q,\delta} \frac{\|\xi\|_{L_q}}{\sqrt{N}}.
\end{equation*}
\endproof

\vskip0.4cm

\noindent{\bf The trace-norm.}

Recall that $B_1$ is the unit ball of the trace norm, that $B_2$ is the unit ball of the Hilbert-Schmidt norm, and that the canonical Gaussian vector here is the Gaussian matrix $G=(g_{ij})$.  Since the operator norm is the dual to the trace norm,
$$
\ell_*(B_1) = \E \sigma_1(G) \lesssim \sqrt{\max\{m,T\}},
$$
and clearly,
$$
\ell_*(B_2) = \E \norm{G}_2 \lesssim \sqrt{mT}.
$$
Thus,
\begin{equation*}
\ell_*(\rho B_\Psi \cap r B_2)=\ell_*(\rho B_1 \cap r B_2)\leq \min\big\{\rho \ell_*(B_1), r \ell_*(B_2)\big\}\lesssim \min\{\rho \sqrt{\max\{m,T\}}, r \sqrt{mT}\}.
\end{equation*}
Therefore,
\begin{equation*}
r_Q^2(\rho) \lesssim_L \left\{
\begin{array}{cc}
0 & \mbox{ if } N \gtrsim_L mT\\
\rho^2\frac{\max\{m,T\}}{N} & \mbox{ otherwise,}
\end{array}
\right.
\end{equation*}
and
\begin{equation*}
r_M^2(\rho) \lesssim_{L,q,\delta} \left\{
\begin{array}{cc}
\norm{\xi}_{L_q}^2\frac{mT}{N} & \mbox{ if } \rho^2 N \gtrsim_{L,q,\delta} \norm{\xi}_{L_q}^2 mT\left(\min{\{m,T\}}\right)^2\\
\\
\rho \|\xi\|_{L_q} \sqrt{\frac{\max\{m,T\}}{N}} & \mbox{ otherwise.}
\end{array}
\right.
\end{equation*}

\noindent{\bf Proof of Theorem \ref{thm:intro-TRACE-est}.}
It is straightforward to verify that if $N \gtrsim_{L,q,\delta} s \max\{m,T\}$ then $s\lesssim (\rho/r(\rho))^2$ when
\begin{equation*}
\rho \gtrsim_{L,q,\delta} \|\xi\|_{L_q} s \sqrt{\frac{\max\{m,T\}}{N}}
\end{equation*} as required in Lemma \ref{lemma:norming-for-trace-norm}.  Moreover, if there is some $V\in\R^{m\times T}$ for which $\norm{V-A^*}_1\lesssim \rho$ and ${\rm rank}(V)\leq s$, it follows that $\Delta(\rho)\geq 4\rho/5$. Setting
\begin{equation*}
\lambda\sim \frac{r^2(\rho^*)}{\rho^*}\sim_{L,q,\delta} \|\xi\|_{L_q} \sqrt{\frac{\max\{m, T\}}{N}},
\end{equation*}
Theorem \ref{lemma:basic-combining-loss-and-reg} yields the bounds on $\|\hat{A}-A^*\|_1$ and $\|\hat{A}-A^*\|_2$. The bounds on the Schatten norms $\|\hat{A}-A^*\|_p$ for $1\leq p\leq 2$ hold because $\norm{A}_p\leq \norm{A}_1^{-1+2/p}\norm{A}_2^{2-2/p}$.
\endproof

\section{Concluding Remarks} \label{sec:concluding-remarks}
As noted earlier, the method we present may be implemented in classical regularization problems as well, leading to an error rate that depends on $\Psi(f^*)$ -- by applying the trivial bound on $\Delta(\rho)$ when $\rho \sim \Psi(f^*)$.

The key issue in classical regularization schemes is the price that one has to pay for not knowing $\Psi(f^*)$ in advance. Indeed, given information on $\Psi(f^*)$, one may use a learning procedure taking values in $\{f \in F : \Psi(f) \leq \Psi(f^*)\}$ such as Empirical Risk Minimization. This approach would result in an error rate of $r(c\Psi(f^*))$, and the hope is that the error rate of the regularized procedure is close to that -- without having prior knowledge on $\Psi(f^*)$. Surprisingly, as we show in \cite{LM_reg_comp}, that is indeed the case.

The problem with applying Theorem \ref{lemma:basic-combining-loss-and-reg} to the classical setup is the choice of $\lambda$. One has no information on $\Psi(f^*)$, and thus setting $\lambda \sim r^2(\rho)/\rho$ for $\rho \sim \Psi(f^*)$ is clearly impossible.

A first attempt of bypassing this obstacle is Remark \ref{rem:large-rho}: if $\rho \gtrsim \Psi(f^*)$, there is no upper constraint on the choice of $\lambda$. Thus, one may consider $\lambda \sim \sup_{\rho > 0} \frac{r^2(\rho)}{\rho}$, which suits any $\rho>0$. Unfortunately, that choice will not do, because in many important examples the supremum happens to be infinite. Instead, one may opt for the lower constraint on $\lambda$ and select
\begin{equation}\label{eq:lambda-via-gamma}
\lambda \sim \sup_{\rho>0} \frac{\gamma_{\cal O}(\rho)}{\rho},
\end{equation}
which is also a legitimate choice for any $\rho$, and is always finite.

We will show in \cite{LM_reg_comp} that the choice in \eqref{eq:lambda-via-gamma} leads to optimal bounds in many interesting examples -- thanks to the first part of Theorem \ref{lemma:basic-combining-loss-and-reg}.

\vskip0.4cm

An essential component in the analysis of regularization problems is bounding $r(\rho)$, and we only considered the subgaussian case and completely ignored the question of the probability estimate. In that sense, the method we presented falls short of being completely satisfactory.

Addressing both these issues requires sharp upper estimates on empirical and multiplier processes, preferably in terms of some natural geometric feature of the underlying class. Unfortunately, this is a notoriously difficult problem. Indeed, the final component in the chaining-based analysis used to study empirical and multiplier processes is to translate a metric complexity parameter (e.g., Talagrand's $\gamma$-functionals) to a geometric one (for example, the mean-width of the set). Such estimates are known almost exclusively in the Gaussian case -- which is, in a nutshell, Talagrand's {\it Majorizing Measures theory} \cite{MR3184689}.

The chaining process in \cite{shahar_multi_pro} is based on a more sensitive metric parameter than the standard Gaussian one. This leads to satisfactory results for other choices of random vectors that are not necessarily subgaussian, for example, unconditional log-concave random vectors. Still, it is far from a complete theory -- as a general version of the Majorizing Measures Theorem is not known.

Another relevant fact is from \cite{Men-multi-weak}. It turns out that if $V$ is a class of linear functionals on $\R^d$ that satisfies a relatively minor symmetry property, and $X$ is an isotropic random vector for which
\begin{equation}\label{eq:sub-gauss-moments}
\sup_{t \in S^{d-1}} \|\inr{X,t}\|_{L_p} \leq L\sqrt{p} \ \ {\rm for} \ \ 2 \leq p \lesssim \log d,
\end{equation}
then the empirical and multiplier processes indexed by $V$ behave as if $X$ were a subgaussian vector. In other words, for such ``symmetric" problems it suffices to have a subgaussian moment growth up to $p \sim \log d$ to ensure a subgaussian behaviour.

This fact is useful because all the indexing sets considered here (and in many other sparsity-based regularization procedures as well) satisfy the required symmetry property.

\vskip0.4cm

Finally, a word about the probability estimate in Theorem \ref{thm:multi-main}. The actual result from \cite{shahar_multi_pro} leads to a probability estimate governed by two factors: the $L_q$ space to which $\xi$ belongs and the `effective dimension' of the class. For a class of linear functionals on $\R^d$ and an isotropic vector $X$, this effective dimension is
$$
D(V)=\left(\frac{\ell^*(V)}{d_2(V)}\right)^2,
$$
where $\ell_*(V)=\E \sup_{v \in V} |\inr{G,v}|$ and $d_2(V)=\sup_{v \in V} \|v\|_{\ell_2^d}$.

One may show that with probability at least
$$
1-c_1w^{-q} N^{-((q/2)-1)}\log^{q} N-2\exp(-c_2u^2 D(V)),
$$
\begin{equation} \label{eq:multi-subgaussian-full}
\sup_{v \in V} \left|\frac{1}{\sqrt{N}}\sum_{i=1}^N \left(\xi_i \inr{V,X_i} - \E \xi \inr{X,v}\right) \right| \lesssim L wu \|\xi\|_{L_q}\ell_*(V).
\end{equation}
If $\xi$ has better tail behaviour, the probability estimate improves; for example, if $\xi$ is subgaussian then \eqref{eq:multi-subgaussian-full} holds with probability at least $1-2\exp(-cw^2N)-2\exp(-cu^2 D(V))$.

The obvious complication is that one has to obtain a {\it lower bound} on the effective dimension $D(V)$. And while it is clear that $D(v) \gtrsim 1$, in many cases (including our three examples) a much better bound is true.

Let us mention that the effective dimension is perhaps the most important parameter in Asymptotic Geometric Analysis. Milman's version of Dvoretzky's Theorem (see, e.g., \cite{MR3331351}) shows that $D(V)$ captures the largest dimension of a Euclidean structure hiding in $V$. In fact, this geometric observation exhibits why that part of the probability estimate in \eqref{eq:multi-subgaussian-full} cannot be improved.


\begin{footnotesize}
\bibliographystyle{plain}
\bibliography{biblio}
\end{footnotesize}

\section{Supplementary material: non-isotropic design} 
\label{sec:supplementary_material_non_isotropic_design}
An inspection of Theorem~\ref{lemma:basic-combining-loss-and-reg} reveals no mention of an isotropicity assumption. There is no choice of a Euclidean structure, and in fact, the statement itself is not even finite dimensional. All that isotropicity has been used for was to bound the ``complexity function'' $ r(\cdot)$  and the ``sparsity function'' $\Delta(\cdot)$ in the three applications --- the LASSO (in Theorem~\ref{thm:intro-LASSO-est}), SLOPE (in Theorem~\ref{thm:intro-SLOPE-est}) and the trace norm regularization (in Theorem~\ref{thm:intro-TRACE-est}).
We may apply Theorem~\ref{lemma:basic-combining-loss-and-reg} to situations that do not involve an isotropic vector and here we give an example of how this may be done.

To simplify our presentation we will only consider $\ell_1$ and SLOPE regularization, which may both be written as 
\begin{equation*}
 \Psi(t) = \sum_{j=1}^d \beta_j t_j^\sharp,
 \end{equation*} 
 where $\beta_1\geq \cdots\geq \beta_d>0$ and $t_1^\sharp\geq \cdots\geq t_d^\sharp\geq0$ is the nondecreasing rearrangement of $(|t_j|)$. As mentioned previously, the LASSO case is recovered for $\beta_1=\cdots=\beta_d=1$ and the SLOPE norm is obtained for $\beta_j = C \sqrt{\log(ed/j)}$ for some constant $C$. We also denote by $B_\Psi$ (resp. $S_\Psi$) the unit ball (resp. sphere) associated with the $\Psi$-norm.

Let $\Sigma\in\R^{d \times d}$ be the covariance matrix of $X$ and set $D=\{x\in\R^d : \norm{\Sigma^{1/2}x}_2\leq 1\}$ to be the corresponding ellipsoid. Naturally, if $X$ is not isotropic than $\Sigma$ is not the identity matrix.

In order to apply Theorem~\ref{lemma:basic-combining-loss-and-reg}, we need to bound from above the expectation of the supremum of the Gaussian process indexed by $\rho B_\Psi\cap r D$:
\begin{equation}\label{eq:gauss_width_non_iso}
\ell_*\left(\rho B_\Psi\cap r D\right) = \E \sup_{w\in \rho B_\Psi\cap r D} \inr{\Sigma^{1/2}G, w}
\end{equation}
where $G$ is a standard Gaussian vector in $\R^d$. 

We also need to solve the ``sparsity equation''---that is, find $\rho^*>0$ for which $\Delta(\rho^*)\geq 4\rho^*/5$ where, for every $\rho>0$,
\begin{equation*}
\Delta(\rho)  = \inf_{h\in\rho S_\Psi\cap r D}\sup_{g\in\Gamma_{t^*}(\rho)}\inr{h, g}
\end{equation*}and  $\Gamma_{t^*}(\rho)$ is the collection of all subgradients of $\Psi$ of vectors in  $t^*+(\rho/20)B_\Psi$.

\vskip0.4cm
We will show that the same results that have been obtained for the LASSO and SLOPE in Theorem~\ref{thm:intro-LASSO-est} and Theorem~\ref{thm:intro-SLOPE-est} actually hold under the following assumption.

\begin{Assumption}\label{assum:covariance_mat}
Let $j \in \{1,...,d\}$ and denote by $\Sigma^{1/2}_{j\bullet}$ the $j$-th row of $\Sigma^{1/2}$. Let $s\in\{1, \ldots, d\}$ and set $\cB_s=\sum_{j=1}^s \beta_j/\sqrt{j}$.
\begin{enumerate}
    \item There exists $\sigma>0$ such that for all $j\in\{1, \ldots, d\}$,  $\norm{\Sigma^{1/2}_{j\bullet}}_2\leq \sigma$.
    \item For all $x\in (20\cB_s S_\Psi)\cap D$, $2\norm{\Sigma^{1/2}x}_2\geq \sup_{|J|\leq s}\norm{x_J}_2$.
\end{enumerate}
\end{Assumption}

\vskip0.4cm

\noindent{\bf SLOPE.}

We first control the Gaussian mean width in \eqref{eq:gauss_width_non_iso} when $\Psi(\cdot)$ is the SLOPE norm.

\begin{Lemma}\label{lem:rearrangement_gaussian}
Set $\beta_j=C\sqrt{\log(ed/j)}$ and let $\Sigma$ be a $d\times d$ symmetric nonnegative matrix for which $\max_j \norm{\Sigma^{1/2}_{j\bullet}}_2\leq \sigma$. If $D=\{x\in\R^d : \norm{\Sigma^{1/2}x}_2\leq 1\}$ then
\begin{equation*}
\E \sup_{w\in \rho B_\Psi\cap r D} \inr{\Sigma^{1/2}G, w}\leq \min\left\{\frac{\rho}{C}\left(\frac{3\sqrt{6}\sigma}{8}+r\sqrt{\frac{\pi}{2}}\right), r \sqrt{d}\right\}
\end{equation*}
\end{Lemma}
\proof
Note that
\begin{equation*}
\E \sup_{w\in \rho B_\Psi\cap r D} \inr{\Sigma^{1/2}G, w}\leq r \E \sup_{w\in  D} \inr{G, \Sigma^{1/2}w} = r\E \norm{G}_2\leq r \left(\E \norm{G}_2^2\right)^{1/2} \leq r \sqrt{d}.
\end{equation*}

Next, let $H:\R^d\to\R$ be defined by $H(u) = \sup\left(\inr{\Sigma^{1/2} u, w}: w \in \rho B_\Psi\cap r D\right)$ and recall that $G$ is the standard Gaussian vector in $\R^d$. It is straightforward to verify that
 \begin{equation*}
 H(G)\leq \sup_{w\in \rho B_\Psi}\inr{\Sigma^{1/2}G, w}\leq \frac{\rho}{C}\max_{1\leq j\leq d}\frac{\xi_j^\sharp}{ \sqrt{\log(ed/j)}}
 \end{equation*}
where we set $(\xi_j)_{j=1}^d = \Sigma^{1/2}G$ and $(\xi_j^\sharp)_{j=1}^d$ is the non-increasing rearrangement of $(|\xi_j|)_{j=1}^d$. Observe that for $u,v\in\R^d$,
 \begin{equation*}
 |H(u) - H(v)|\leq \sup_{w \in \rho B_\Psi\cap r D}|\inr{\Sigma^{1/2}(u-v),w}| \leq  r \sup_{w\in B_2^d}|\inr{u-v, w}|= r \norm{u-v}_2,
 \end{equation*}
implying that $H$ is a Lipschitz function with constant $r$; thus, it follows from p.~21 in Chapter~1 of \cite{MR1849347} that 
\begin{equation}\label{eq:mean_median}
\E H(G)\leq {\rm Med}(H(G)) + r\sqrt{\frac{\pi}{2}},
\end{equation}
where ${\rm Med}(H(G))$ is the median of $H(G)$.

Hence, to obtain the claimed bound on $\E \sup_{w\in \rho B_\Psi}\inr{\Sigma^{1/2}G, w}$ it suffices to establish a suitable upper estimate on the median of $\max_{1\leq j\leq d}\frac{\xi_j^\sharp}{ \sqrt{\log(ed/j)}}$. With that in mind, let $\xi_1,\ldots, \xi_N$ be mean-zero Gaussian variables and assume that for every $j=1, \ldots, d$,
\begin{equation}\label{eq:subgauss_var}
\E \exp(\xi_j^2/L^2)\leq e
\end{equation}
for some $L>0$. Note that in our case, $\xi_1,\ldots, \xi_N$ satisfying \eqref{eq:subgauss_var} for $L = 3 \sigma/8$. 

By Jensen's inequality,
\begin{equation*}
\E \exp\left(\frac{1}{j}\sum_{k=1}^j \frac{(\xi_k^\sharp)^2}{L^2}\right)\leq \frac{1}{j}\sum_{k=1}^j\E \exp\left( \frac{(\xi_k^\sharp)^2}{L^2}\right)\leq \frac{1}{j}\sum_{k=1}^d\E \exp\left( \frac{\xi_k^2}{L^2}\right)\leq \frac{ed}{j};
\end{equation*}
hence,
\begin{equation}\label{eq:ordered_subgauss_1}
Pr\left(\frac{1}{j}\sum_{k=1}^j \frac{(\xi_k^\sharp)^2}{L^2}\geq 2 \log\left(ed/j\right)\right)\leq \exp\left(-\log\left(ed/j\right)\right)=\frac{j}{ed}.
\end{equation}

Let $q\geq0$ be the integer that satisfies $2^q\leq d < 2^{q+1}$. It follows from \eqref{eq:ordered_subgauss_1} that with probability at least
\begin{equation*}
1-\sum_{\ell=0}^{q-1}\frac{2^\ell}{ed} = 1-\frac{2^{q}-1}{ed}>\frac{1}{2},
\end{equation*}for every $\ell=0,\cdots, q-1$,
\begin{equation*}
 (\xi_{2^\ell}^\sharp)^2\leq \frac{1}{2^{\ell}}\sum_{j=1}^{2^\ell}(\xi_{j}^\sharp)^2\leq 2 L^2\log(ed/2^\ell).
 \end{equation*} 
 Moreover, for $2^{\ell}\leq j < 2^{\ell+1}$, we have $\xi^\sharp_j\leq \xi_{2^\ell}^\sharp$ and $\log(ed/2^\ell)\leq 2 \log(ed/j)$; also for $2^q\leq j\leq d$, we have $\xi_j^\sharp\leq \xi_{2^{q-1}}^\sharp$ and $\log(ed/2^\ell)\leq 3 \log(ed/j)$. Therefore,
 \begin{equation}\label{eq:mediane}
 Pr \left(\max_{1\leq j\leq d}\frac{\xi_j^\sharp}{ \sqrt{\log(ed/j)}}\leq \sqrt{6}L\right)>\frac{1}{2},
 \end{equation}
proving the requested bound on ${\rm Med}(H(G))$.
\endproof

Observe that up to constant $\sigma$, we actually recover the same result as in \eqref{eq:gauss_mean_slope_iso}; therefore, one may choose the same ``complexity function'' $r(\cdot)$ as in the proof of Theorem~\ref{thm:intro-SLOPE-est}. 

\vskip0.4cm
Let us turn to a lower bound on the ``sparsity function''.
\begin{Lemma}\label{lem:sparsity_function_slope}
There exists an absolute constant $0<c<80$ for which the following holds. Let $s\in\{1, \ldots, d\}$ and set $\cB_s = \sum_{j\leq s}\beta_j/\sqrt{j}$.  Assume that for every $x\in (80 \cB_s S_\psi) \cap D$ one has $\norm{\Sigma^{1/2}x}_2\geq (1/2)\sup_{|J|\leq s}\norm{x_J}_2$. Let $\rho>0$ and assume further that there is a $s$-sparse vector in $t^*+(\rho/20) B_\Psi$. If $80\cB_s\leq \rho/r(\rho)$ then
\begin{equation*}
\Delta(\rho)  = \inf_{h\in\rho S_\Psi\cap r(\rho) D}\sup_{g\in\Gamma_{t^*}(\rho)}\inr{h, g}\geq \frac{4\rho}{5}.
\end{equation*}
\end{Lemma}

\proof Let $h\in\rho S_\Psi\cap r(\rho) D$ and denote by $(h_j^\sharp)$ the non-increasing rearrangement of $(|h_j|)$. It follows from the proof of Lemma~\ref{lemma:delta-for-slope} that
\begin{equation*}
\sup_{g\in\Gamma_{t^*}(\rho)}\inr{h, g}\geq \frac{17\rho}{20}-2 \sum_{j\leq s}\beta_j h_j^\sharp.
\end{equation*}
Let $h^{\sharp,s}$ be the $s$-sparse vector with coordinates given by $h_j^\sharp$ for $1\leq j\leq s$ and $0$ otherwise. We have
\begin{equation*}
\frac{h}{\rho}\in S_\Psi\cap \left(\frac{r(\rho)}{\rho}\right) D \subset S_\Psi \cap \left(\frac{1}{80 \cB_s}\right) D
\end{equation*}
implying that $2\norm{\Sigma^{1/2} h}_2\geq \norm{h^{\sharp, s}}_2$. Furthermore, since $h_j^\sharp\leq \norm{h^\sharp}_2/\sqrt{j}$ for every $1\leq j\leq s$, we have
\begin{equation*}
\sum_{j\leq s}\beta_j h_j^\sharp\leq \norm{h^{\sharp,s}}_2\cB_s\leq 2\cB_s \norm{\Sigma^{1/2}h}_2\leq 2\cB_s r(\rho).
\end{equation*}Hence, if $\rho\geq80 r(\rho)\cB_s$ then $\Delta(\rho)\geq 4 \rho/5$.
\endproof

We thus recover the same condition as in Lemma~\ref{lemma:delta-for-slope}, implying that  Theorem~\ref{thm:intro-SLOPE-est} actually holds under the weaker Assumption~\ref{assum:covariance_mat}: let $X$ be an $L$-subgaussian random vector whose covariance matrix satisfies Assumption~\ref{assum:covariance_mat}. The SLOPE estimator with regularization parameter $\lambda\sim \norm{\xi}_{L_q}/\sqrt{N}$ satisfies, with probability at least $1-\delta-\exp(-c_0NL^8)$,
\begin{equation*}
\Psi(\hat t-t^*)\leq c_3\|\xi\|_{L_q} \frac{s}{\sqrt{N}} \log\Big(\frac{ed}{s}\Big) \ \ \mbox{ and } \ \ \norm{\Sigma^{1/2}(\hat t-t^*)}_2^2 \leq c_3 \|\xi\|_{L_q}^2 \frac{s }{N}\log\Big(\frac{ed}{s}\Big)
\end{equation*} 
when $N\geq c_4 s \log(ed/s)$ and when there is a s-sparse vector close enough to $t^*$.

\vskip0.4cm

\noindent{\bf The LASSO.}

Here, for every $1 \leq j \leq d$, $\beta_j=1$;  $B_\Psi = B_1^d$; and $\cB_s=\sum_{j=1}^s 1/\sqrt{j}\leq 2 \sqrt{s}$.

\begin{Lemma}\label{lem:rearrangement_gaussian_lasso}
Let $\Sigma$ be a $d\times d$ symmetric nonnegative matrix for which $\max_j \norm{\Sigma^{1/2}_{j\bullet}}_2\leq \sigma$. If $D=\{x\in\R^d : \norm{\Sigma^{1/2}x}_2\leq 1\}$ then every $\rho>0$ and $r>0$,
\begin{equation*}
\E \sup_{w\in \rho B_1^d\cap r D} \inr{\Sigma^{1/2}G, w}\leq \min\left\{r\sqrt{d}, \rho\sigma \sqrt{\log(ed)}\right\}.
\end{equation*}
\end{Lemma}
\proof
Note that
\begin{equation*}
\E \sup_{w\in \rho B_1^d\cap r D} \inr{\Sigma^{1/2}G, w}\leq r \E \sup_{w\in D} \inr{G, \Sigma^{1/2}w} = r\E \norm{G}_2\leq r \left(\E \norm{G}_2^2\right)^{1/2} \leq r \sqrt{d}.
\end{equation*}
Next, set $(\xi_j)_{j=1}^d = \Sigma^{1/2}G$ and let $(\xi_j^\sharp)_{j=1}^d$ be the non-increasing rearrangement of $(|\xi_j|)$. Therefore, $\xi_1,\ldots, \xi_d$ are mean-zero Gaussian variables and satisfy $\E \exp(\xi_j^2/L^2)\leq e$ for $L = 3 \sigma/8$. It is evident that 
\begin{align}\label{eq:rearrangement}
\nonumber &\E\left(\frac{(\xi_1^\sharp)^2}{L^2}\right) \leq \log \left(\E\exp\left(\frac{(\xi_1^\sharp)^2}{L^2}\right)\right)\leq \log \left(\sum_{j=1}^d \E\exp\left(\frac{\xi_j^2}{L^2}\right)\right)\leq \log\left(ed\right),
\end{align}
and therefore,
\begin{align*}
\E \sup_{w\in \rho B_1^d\cap r D} \inr{\Sigma^{1/2}G, w}\leq \E \sup_{w\in \rho B_1^d } \inr{\Sigma^{1/2}G, w}\leq \rho \E\left((\xi_1^\sharp)^2\right)^{1/2}\leq \frac{3\rho\sigma L}{8} \sqrt{\log\left(ed\right)}.
\end{align*}
\endproof

Lemma~\ref{lem:rearrangement_gaussian_lasso} leads to a slightly different result than in the isotropic case (Lemma~\ref{lemma:monotone-for-slope}), and as a consequence, $r(\cdot)$ has to be slightly modified. A straightforward computation shows that
\begin{equation*}
r_M^2(\rho) \lesssim_{L,q,\delta} \min\left(\frac{\norm{\xi}_{L_q} d }{N} ,  \rho \sigma \norm{\xi}_{L_q} \sqrt{\frac{\log(ed)}{N}}\right)
\end{equation*}
and
\begin{equation*}
r_Q^2(\rho) \lesssim_L  \left\{
\begin{array}{cc}
  0 & \mbox{ if }  N \gtrsim_L  d \\
\frac{\rho^2\sigma^2}{N}\log\Big(\frac{c(L)d}{N}\Big) & \mbox{ otherwise}.
\end{array}
\right.
\end{equation*}
and still $r(\rho) = \max\{r_M(\rho), r_Q(\rho)\}$.

Finally, let us prove the sparsity condition.

\begin{Lemma}\label{lem:sparsity_function_lasso}
There exists an absolute constant $0<c<80$ for which the following holds. Let $s\in\{1, \ldots, d\}$ and set $\cB_s = \sum_{j\leq s}1/\sqrt{j}$.  Assume that for every $x\in \left(20\cB_s S(\ell_1^d)\right)\cap D$ one has $\norm{\Sigma^{1/2}x}_2\geq (1/2)\sup_{|J|\leq s}\norm{x_J}_2$. Let $\rho>0$ and assume further that there is a $s$-sparse vector in $t^*+(\rho/20) B_1^d$. If $20\cB_s\leq \rho/r(\rho)$ then
\begin{equation*}
\Delta(\rho)  = \inf_{h\in\rho S(\ell_1^d)\cap r(\rho) D}\sup_{g\in\Gamma_{t^*}(\rho)}\inr{h, g}\geq \frac{4\rho}{5}.
\end{equation*}
\end{Lemma}

\proof Let $h\in\rho S(\ell_1^d)\cap r(\rho) D$ and denote by $(h_j^\sharp)$ the non-increasing rearrangement of $(|h_j|)$. It follows from the proof of Lemma~\ref{lemma:norming-for-LASSO} that
\begin{equation*}
\sup_{g\in\Gamma_{t^*}(\rho)}\inr{h, g}\geq \rho-2 \sum_{j\leq s} h_j^\sharp.
\end{equation*}
Let $h^{\sharp,s}$ be the $s$-sparse vector with coordinates given by $h_j^\sharp$ for $1\leq j\leq s$ and $0$ otherwise. Observe that
\begin{equation*}
\frac{h}{\rho}\in S(\ell_1^d)\cap \left(\frac{r(\rho)}{\rho}\right)D\subset S(\ell_1^d)\cap \left(\frac{1}{40 \cB_s}\right)D,
\end{equation*}
and therefore $2\norm{\Sigma^{1/2}h}_2\geq \norm{h^{\sharp, s}}_2$.

It follows that 
\begin{equation*}
\sum_{j\leq s}h_j^\sharp\leq \sqrt{s}\norm{h^\sharp}_2\leq \cB_s\norm{h^\sharp}_2\leq 2\cB_s \norm{\Sigma^{1/2}h}_2\leq 2\cB_s r(\rho),
\end{equation*}
and in particular, if $\rho\geq20 r(\rho)\cB_s$ then $\Delta(\rho)\geq 4 \rho/5$.
\endproof

Using the estimate on $r(\cdot)$ and  Lemma~\ref{lem:sparsity_function_lasso}, it is evident that when $N\gtrsim_{L, q, \delta} s \sigma^2 \log(ed)$, one has $\Delta(\rho^*)\geq 4 \rho^*/5$ for
\begin{equation*}
\rho^*\sim_{L, q, \delta} \norm{\xi}_{L_q} s \sqrt{\frac{\log(ed)}{N}}
\end{equation*}
and if there is a $s$-sparse vector in $t^*+(\rho^*/20)B_1^d$. 

Finally, one may choose the regularization parameter by setting
\begin{equation*}
\lambda\sim \frac{r^2(\rho^*)}{\rho^*} \sim_{L, q, \delta} \norm{\xi}_{L_q} s \sqrt{\frac{\log(ed)}{N}}.
\end{equation*}
It follows that if $X$ is an $L$-subgaussian random vector that satisfies Assumption~\ref{assum:covariance_mat} then with probability larger than $1-\delta-2 \exp(-c_0NL^8)$,
\begin{equation*}
\norm{\hat t - t^*}_1\leq \rho^* = c_1(\delta)\norm{\xi}_{L_q} s \sqrt{\frac{\log(ed)}{N}} \mbox{ and } \norm{\Sigma^{1/2}(\hat t - t^*)}_2\leq r(\rho^*)=c_2(L,\delta)\norm{\xi}_{L_q} \sqrt{\frac{s\log(ed)}{N}}.
\end{equation*}


\end{document}